\documentclass{amsart}
\usepackage[dvips]{epsfig}
\usepackage[dvips]{graphics}
\usepackage[dvips]{color}
\usepackage{amsfonts}
\usepackage{amsmath}
\usepackage{amssymb}
\usepackage{amsthm}
\usepackage{graphics}
\newtheorem{thm}{Theorem}[section]

 \newtheorem{prop}[thm]{Proposition}
 \theoremstyle{definition}

 \newtheorem{rem}{Remark}
 \numberwithin{equation}{section}

\begin{document}
\title[Positive Dehn Twist Expression for a $\mathbb{Z}_3$ action]
 {Positive Dehn Twist Expression for a $\mathbb{Z}_3$ action on $\Sigma_g$}

\author{ HISAAKI ENDO AND YUSUF Z. GURTAS }

 \address{Department of Mathematics, Graduate School of Science, Osaka University,
Toyonaka, Osaka 560-0043, Japan}

\email{endo@math.sci.osaka-u.ac.jp}

 \address{Department of Mathematics and Computer Science, St. Louis University, MO, USA}

 \email{ygurtas@slu.edu}

\thanks{}




\dedicatory{}

\commby{}

\begin{abstract}
A positive Dehn twist product for a $\mathbb{Z}_3$ action 
 on the 2-dimensional closed, compact, oriented
surface $\Sigma_g$ is presented. The homeomorphism invariants of the
resulting   symplectic 4-manifolds are computed.

\end{abstract}

\maketitle

\section*{Introduction}
 This article attempts to answer a question raised by Feng Luo in
\cite{Lu3} which asks for a  Dehn twist expression for the generator
of a
 $\mathbb{Z}_3$ action with $g+2$
fixed points on the 2-dimensional closed, compact, oriented surface
$\Sigma_g$. By the work of Nielsen, there is only one such action on
$\Sigma_g$, \cite{Lu3}. \\
\indent In Section \ref{construction} we build a closed genus $g$
surface $\Sigma_g$ using $g$ copies of tori with boundary as
building blocks in order to realize that action on $\Sigma_g$. We
simply take an order three element from the mapping class group
$\mathcal{M}_1$ of torus and juxtapose its Dehn twist expression in
$\mathcal{M}_g$, considering torus with boundary as a subsurface of
$\Sigma_g$ and taking the orientation into consideration in the
gluing process. We start with a torus with one boundary component
oriented positively. Then glue a torus with two boundary components
oriented negatively to it. Then keep adding more tori with boundary
with alternating orientations and finally cap it off with a torus
with one boundary component oppositely oriented as the previous
copy. We aim at a Dehn twist product for the generator of the
$\mathbb{Z}_3$ action on $\Sigma_g$ that uses only positive
exponents  in order to make sure that the 4-manifold it defines as
Lefschetz fibration carries symplectic structure. This becomes a
challenge because the negatively oriented bounded tori introduce
into the expression many elements with negative exponents and there
are still some negative powers to be eliminated from the expression
for genus $g>6$. Therefore this work is still in progress.\\

\indent In Section \ref{low genus} we show explicitly how to obtain
a positive Dehn twist product for the generator of the
$\mathbb{Z}_3$ action on $\Sigma_g,g\leq6.$ What seems to be working
for low genus doesn't generalize to higher genus easily and
the construction evolves rather ad hoc, at least partially. \\
\indent In Section \ref{applications} we compute the Euler
characteristic and signatures of the 4-manifolds given by the words
that are obtained in  Section \ref{low genus}. The method introduced
by the first author and S.Nagami is used for  signature
computations, \cite{ES}.

\section{Review of Relations in $\mathcal{M}_1,\mathcal{M}_1^1$ and
$\mathcal{M}_1^2$ }

\begin{figure}[htbp]
     \centering  \leavevmode
     \psfig{file=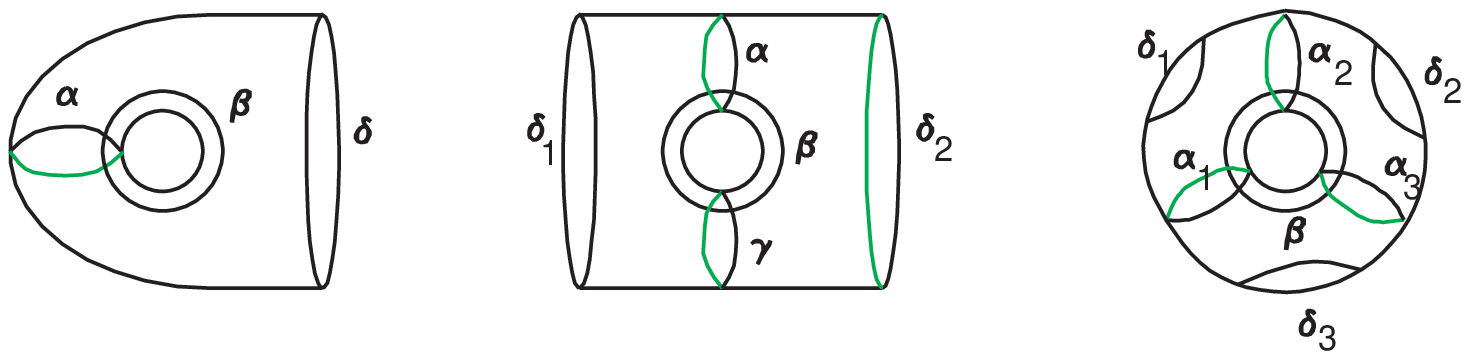,width=5.80in,clip=}
     \caption{} \label{figure-torus}
 \end{figure}

The mapping class group $\mathcal{M}_1$ of  torus is generated by
Dehn twists about the cycles $\alpha$ and $\beta$, Figure
\ref{figure-torus}, subject to the relations
\begin{eqnarray} \label{relations-torus}
\nonumber \alpha \beta \alpha &=&\beta \alpha \beta \\
\left( \alpha \beta \right) ^{6}&=&1.
\end{eqnarray}
Here, by abuse of notation, we use $\alpha$ and $ \beta  $ to mean
Dehn twists about them for simplicity. The first relation is called \emph{braid relation} and it
exists between every pair of curves that intersect transversely.\\
\indent Torus with one, two, and three boundary components are
subject to the relations

\begin{equation} \label{basic relations}
\left( \alpha \beta \right) ^{6}=\delta \ \  \text{and}\ \ \left(
\beta \alpha \beta \gamma  \right) ^{3}= \left(  \alpha \beta \gamma
\right) ^{4}=\delta _{1}\delta_{2}, \ \ \text{and}\ \ \left(
\alpha_1 \alpha_2 \alpha_3 \beta \right) ^{3}=\delta
_{1}\delta_{2}\delta_{3}
\end{equation}
 respectively. The last one is also called \emph{star relation}, \cite{Ge}. \\
 \vspace{-.15in}
 The basic idea that is used in this paper is to
glue several copies of torus with two boundary components together
and cap  the resulting bounded surface off with two copies of torus
with one boundary component, one on each end, to get a closed
surface of genus $g$. We take the word
\begin{equation} \label{word-on-torus-with-one-bd}
\left( \alpha \beta \right) ^{2} \end{equation} on the two end
copies and the word \begin{equation}
\label{word-on-torus-with-two-bd} \beta \alpha \beta \gamma
\end{equation} on each of the remaining copies in between and
juxtapose them with alternating signs to come up with an order three
element in the mapping class group of the resulting closed genus $g$
surface.

\section{Construction of the order three element on  $\Sigma_g$} \label{construction}
 In this section we will construct an order
three element in the mapping class group of  closed genus $g$
surface using the words \eqref{word-on-torus-with-one-bd} and
\eqref{word-on-torus-with-two-bd}
according to their position in the gluing process. First case is when genus $g$ is even. \\

\subsection{Genus g-even} \label{genus-g-even}

\begin{figure}[htbp] \label{figure-genus 1 broken from }
     \centering  \leavevmode
     \psfig{file=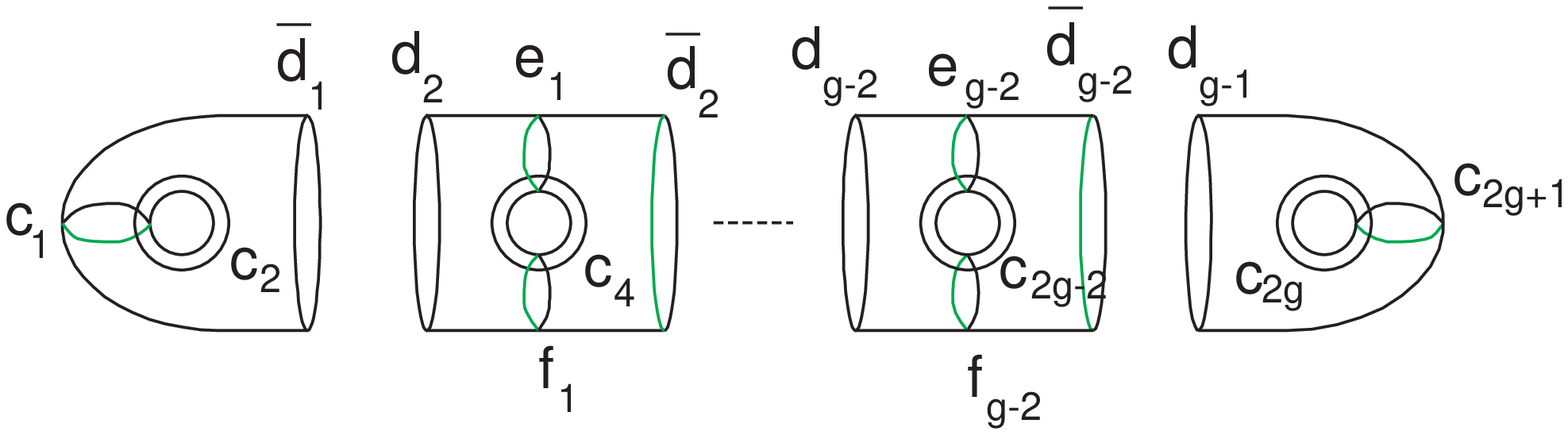,width=5.80in,clip=}
    \caption{}
 \end{figure}

We juxtapose the words of type \eqref{word-on-torus-with-one-bd} and
\eqref{word-on-torus-with-two-bd} on each of the bounded surfaces
above by paying careful attention to the orientation:

 \begin{eqnarray} \label{words-before-gluing-even}
\nonumber \left( c_{1}c_{2}\right) ^{2} \\
\nonumber \left( c_{4}e_{1}c_{4}f_{1}\right) ^{-1} \\
\nonumber c_{6}e_{2}c_{6}f_{2} \\
\nonumber \left( c_{8}e_{3}c_{8}f_{3}\right) ^{-1} \\
 \vdots  \hspace{.5in}  \\
\nonumber \left( c_{2g-4}e_{g-3}c_{2g-4}f_{g-3}\right)^{-1}\\
\nonumber c_{2g-2}e_{g-2}c_{2g-2}f_{g-2} \\
 \nonumber \left( c_{2g+1}c_{2g}\right) ^{-2}
 \end{eqnarray}

Every other surface will be negatively oriented so that we can glue
the boundaries together. Using the chain relation \[\left(
c_{2i+2}e_{i}c_{2i+2}f_{i}\right) ^{3}=d_{i}\bar{d}_{i}\] on torus
with two boundary components and \[\left( c_{2g+1}c_{2g}\right)
^{6}=d_{g-1}\] on torus with one boundary component, the expressions
containing negative exponents in \eqref{words-before-gluing-even}
can be written as
\begin{eqnarray}\label{d's-and-dbar's-even}
 \nonumber\left( c_{4}e_{1}c_{4}f_{1}\right) ^{-1} = d_2^{-1}\left( c_{4}e_{1}c_{4}f_{1}\right) ^{2}\overline{d}_2^{-1}\\
  \left( c_{8}e_{3}c_{8}f_{3}\right) ^{-1} = d_3^{-1}\left( c_{8}e_{3}c_{8}f_{3}\right) ^{2}\overline{d}_3^{-1} \\
\nonumber \vdots  \hspace{.5in}\\
  \nonumber  \left( c_{2g-4}e_{g-3}c_{2g-4}f_{g-3}\right)^{-1}=d_{g-3}^{-1}\left( c_{2g-4}e_{g-3}c_{2g-4}f_{g-3}\right)^{2}\overline{d}_{g-3}^{-1}
\end{eqnarray}
and the last one as
\begin{eqnarray}
  \label{last-d-even} \left( c_{2g+1}c_{2g}\right) ^{-2}= d_{g-1}^{-1}\left( c_{2g+1}c_{2g}\right)
  ^{4}.
\end{eqnarray}

\begin{figure}[htbp]
     \centering  \leavevmode
     \psfig{file=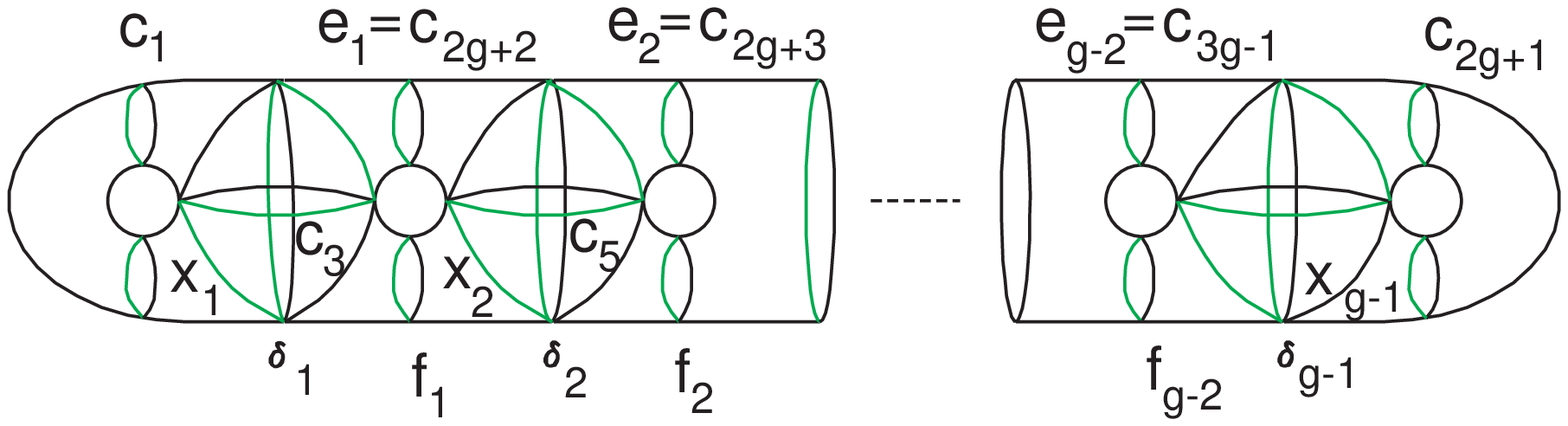,width=5.50in,clip=}
     \caption{} \label{figure-lanterns}
 \end{figure}

 We glue the bounded surfaces together and use the lantern relations
\begin{eqnarray}\label{lanterns-solved-for-delta}
\nonumber \delta _{1}x_{1}c_{3}&=&c_{1}c_{1}e_{1}f_{1} \\
\nonumber \delta _{2}x_{2}c_{5}&=&e_{1}f_{1}e_{2}f_{2} \\
&\vdots&\\
\nonumber \delta _{g-2}x_{g-2}c_{2g-3}&=&e_{g-3}f_{g-3}e_{g-2}f_{g-2} \\
\nonumber \delta
_{g-1}x_{g-1}c_{2g-1}&=&e_{g-2}f_{g-2}c_{2g+1}c_{2g+1}
\end{eqnarray}
 to eliminate the negative exponents of $ \overline{d}_i$ and $d_i$ in \eqref{d's-and-dbar's-even} and \eqref{last-d-even} using the fact that
 $\overline{d}_i= d_{i+1}=\delta_i$ after gluing. Solving \eqref{lanterns-solved-for-delta} for $\delta_i^{-1}$ we get
\begin{eqnarray}\label{lanterns-solved-for-delta-inverse}
\nonumber \delta _{1}^{-1}&=&x_{1}c_{3}c_{1}^{-1}c_{1}^{-1}e_{1}^{-1}f_{1}^{-1} \\
\nonumber \delta _{2}^{-1}&=&x_{2}c_{5}e_{1}^{-1}f_{1}^{-1}e_{2}^{-1}f_{2}^{-1} \\
 &\vdots&\\
\nonumber \delta _{g-2}^{-1}&=&x_{g-2}c_{2g-3}e_{g-3}^{-1}f_{g-3}^{-1}e_{g-2}^{-1}f_{g-2}^{-1} \\
\nonumber
\delta_{g-1}^{-1}&=&x_{g-1}c_{2g-1}e_{g-2}^{-1}f_{g-2}^{-1}c_{2g+1}^{-1}c_{2g+1}^{-1}
\end{eqnarray}\\
 Therefore equations \eqref{d's-and-dbar's-even} and \eqref{last-d-even} become

\begin{eqnarray}\label{d's-eliminated}
 \nonumber \left( c_{4}e_{1}c_{4}f_{1}\right) ^{-1}=x_{1}c_{3}c_{1}^{-1}c_{1}^{-1}e_{1}^{-1}f_{1}^{-1}
 \left( c_{4}e_{1}c_{4}f_{1}\right) ^{2}x_{2}c_{5}e_{1}^{-1}f_{1}^{-1}e_{2}^{-1}f_{2}^{-1}\\
\nonumber  \left( c_{8}e_{3}c_{8}f_{3}\right) ^{-1}=
  x_{3}c_{7}e_{2}^{-1}f_{2}^{-1}e_{3}^{-1}f_{3}^{-1}\left( c_{8}e_{3}c_{8}f_{3}\right) ^{2}x_{4}c_{9}e_{3}^{-1}f_{3}^{-1}e_{4}^{-1}f_{4}^{-1} \\
 \vdots  \hspace{1.5in}
\end{eqnarray}
\begin{eqnarray*}
  \nonumber \left( c_{2g-4}e_{g-3}c_{2g-4}f_{g-3}\right)^{-1}=
  x_{g-3}c_{2g-5}e_{g-4}^{-1}f_{g-4}^{-1}e_{g-3}^{-1}f_{g-3}^{-1}\left( c_{2g-4}e_{g-3}c_{2g-4}f_{g-3}\right)^{2}
  x_{g-2}c_{2g-3}e_{g-3}^{-1}f_{g-3}^{-1}e_{g-2}^{-1}f_{g-2}^{-1}
\end{eqnarray*}
and
\begin{eqnarray}
  \label{last-d-eliminated} \left( c_{2g+1}c_{2g}\right) ^{-2}= x_{g-1}c_{2g-1}e_{g-2}^{-1}f_{g-2}^{-1}c_{2g+1}^{-1}c_{2g+1}^{-1}\left( c_{2g+1}c_{2g}\right)
  ^{4},
\end{eqnarray}

\bigskip

and \eqref{words-before-gluing-even} becomes
\[ \left( c_{1}c_{2}\right) ^{2} \]
\[ x_{1}c_{3}c_{1}^{-1}c_{1}^{-1}e_{1}^{-1}f_{1}^{-1}
 \left( c_{4}e_{1}c_{4}f_{1}\right)
 ^{2}x_{2}c_{5}e_{1}^{-1}f_{1}^{-1}e_{2}^{-1}f_{2}^{-1}\]
 \[ c_{6}e_{2}c_{6}f_{2} \]
 \[ x_{3}c_{7}e_{2}^{-1}f_{2}^{-1}e_{3}^{-1}f_{3}^{-1}\left( c_{8}e_{3}c_{8}f_{3}\right) ^{2}x_{4}c_{9}e_{3}^{-1}f_{3}^{-1}e_{4}^{-1}f_{4}^{-1}
 \]
\[ c_{10}e_{4}c_{10}f_{4} \]
 \begin{eqnarray} \label{words-put-together-even}
 \vdots
 \end{eqnarray}
 \[ c_{2g-6}e_{g-4}c_{2g-6}f_{g-4} \]
\[ x_{g-3}c_{2g-5}e_{g-4}^{-1}f_{g-4}^{-1}e_{g-3}^{-1}f_{g-3}^{-1}\left(
c_{2g-4}e_{g-3}c_{2g-4}f_{g-3}\right)^{2}
  x_{g-2}c_{2g-3}e_{g-3}^{-1}f_{g-3}^{-1}e_{g-2}^{-1}f_{g-2}^{-1}\]
 \[  c_{2g-2}e_{g-2}c_{2g-2}f_{g-2} \]
\[
x_{g-1}c_{2g-1}e_{g-2}^{-1}f_{g-2}^{-1}c_{2g+1}^{-1}c_{2g+1}^{-1}\left(
c_{2g+1}c_{2g}\right)
  ^{4}.\]

Juxtaposing these words, we obtain

\[
 \left( c_{1}c_{2}\right) ^{2}
x_{1}c_{3}c_{1}^{-1}c_{1}^{-1}e_{1}^{-1}f_{1}^{-1}
 \left( c_{4}e_{1}c_{4}f_{1}\right) ^{2}x_{2}c_{5}e_{1}^{-1}f_{1}^{-1}e_{2}^{-1}f_{2}^{-1}
 c_{6}e_{2}c_{6}f_{2} \]\[
  x_{3}c_{7}e_{2}^{-1}f_{2}^{-1}e_{3}^{-1}f_{3}^{-1}\left( c_{8}e_{3}c_{8}f_{3}\right) ^{2}x_{4}c_{9}e_{3}^{-1}f_{3}^{-1}e_{4}^{-1}f_{4}^{-1}
c_{10}e_{4}c_{10}f_{4} \]\[
  x_{5}c_{11}e_{4}^{-1}f_{4}^{-1}e_{5}^{-1}f_{5}^{-1}\left( c_{12}e_{5}c_{12}f_{5}\right) ^{2}x_{6}c_{13}e_{5}^{-1}f_{5}^{-1}e_{6}^{-1}f_{6}^{-1}
c_{14}e_{6}c_{14}f_{6} \]
\begin{eqnarray}\label{juxtaposed word-even}
\vdots  \hspace{.5in}
\end{eqnarray}
\[c_{2g-6}e_{g-4}c_{2g-6}f_{g-4}
x_{g-3}c_{2g-5}e_{g-4}^{-1}f_{g-4}^{-1}e_{g-3}^{-1}f_{g-3}^{-1}
c_{2g-4}e_{g-3}c_{2g-4}f_{g-3}c_{2g-4}e_{g-3}c_{2g-4}f_{g-3}
    \]\[
 x_{g-2}c_{2g-3}e_{g-3}^{-1}f_{g-3}^{-1}e_{g-2}^{-1}f_{g-2}^{-1} c_{2g-2}e_{g-2}c_{2g-2}f_{g-2}
x_{g-1}c_{2g-1}e_{g-2}^{-1}f_{g-2}^{-1}c_{2g+1}^{-1}c_{2g+1}^{-1}\left(
c_{2g+1}c_{2g}\right)
  ^{4}\]\\

Next, we will eliminate the negative exponents using braid and
commutativity relations only. Let's expand the parenthesis in the
top three lines in \eqref{juxtaposed word-even}:

\[
c_{1}\underline{c_{2}c_{1}c_{2}}
x_{1}c_{3}c_{1}^{-1}c_{1}^{-1}e_{1}^{-1}f_{1}^{-1}
 \underline{c_{4}e_{1}c_{4}}f_{1}\underline{c_{4}e_{1}c_{4}}f_{1}x_{2}c_{5}e_{1}^{-1}f_{1}^{-1}e_{2}^{-1}f_{2}^{-1}
\underline{ c_{6}e_{2}c_{6}}f_{2} \]\[
  x_{3}c_{7}e_{2}^{-1}f_{2}^{-1}e_{3}^{-1}f_{3}^{-1}\underline{c_{8}e_{3}c_{8}}f_{3} \underline{c_{8}e_{3}c_{8}}f_{3}
  x_{4}c_{9}e_{3}^{-1}f_{3}^{-1}e_{4}^{-1}f_{4}^{-1}
\underline{c_{10}e_{4}c_{10}}f_{4} \]\[
  x_{5}c_{11}e_{4}^{-1}f_{4}^{-1}e_{5}^{-1}f_{5}^{-1}\underline{c_{12}e_{5}c_{12}}f_{5} \underline{c_{12}e_{5}c_{12}}f_{5}
  x_{6}c_{13}e_{5}^{-1}f_{5}^{-1}e_{6}^{-1}f_{6}^{-1}
\underline{c_{14}e_{6}c_{14}}f_{6} \] \medskip

and use  braid relation for the underlined triples:

\begin{eqnarray}\label{beginning-of-the-juxtaposed-word} \\
\nonumber c_{1}c_{1}c_{2}\underline{c_{1}}
x_{1}c_{3}\underline{c_{1}}^{-1}c_{1}^{-1}\underline{e_{1}}^{-1}f_{1}^{-1}
 \underline{e_{1}}c_{4}e_{1}f_{1}e_{1}c_{4}\underline{e_{1}f_{1}}x_{2}c_{5}\underline{e_{1}^{-1}f_{1}}^{-1}\underline{e_{2}}^{-1}f_{2}^{-1}
\underline{e_{2}}c_{6} \underline{e_{2}f_{2}}\\
\nonumber
x_{3}c_{7}\underline{e_{2}^{-1}f_{2}}^{-1}\underline{e_{3}}^{-1}f_{3}^{-1}\underline{e_{3}}c_{8}e_{3}f_{3}
e_{3}c_{8}
  \underline{e_{3}f_{3}}x_{4}c_{9}\underline{e_{3}^{-1}f_{3}}^{-1}\underline{e_{4}}^{-1}f_{4}^{-1}
\underline{e_{4}}c_{10}\underline{e_{4}f_{4}}
\end{eqnarray}
\[ x_{5}c_{11}\underline{e_{4}^{-1}f_{4}}^{-1}\underline{e_{5}}^{-1}f_{5}^{-1}\underline{e_{5}}c_{12}e_{5}f_{5}
e_{5}c_{12}
  \underline{e_{5}f_{5}}x_{6}c_{13}\underline{e_{5}^{-1}f_{5}}^{-1}\underline{e_{6}}^{-1}f_{6}^{-1}
\underline{e_{6}}c_{14}e_{6}f_{6}.\] \medskip

Next, cancel the underlined pairs above  using commutativity:

\[
c_{1}c_{1}c_{2} x_{1}c_{3}c_{1}^{-1}f_{1}^{-1}
c_{4}e_{1}f_{1}e_{1}c_{4}x_{2}c_{5}f_{2}^{-1} c_{6}
  x_{3}c_{7}f_{3}^{-1}c_{8}e_{3}f_{3} e_{3}c_{8}
x_{4}c_{9}f_{4}^{-1} c_{10}
  x_{5}c_{11}f_{5}^{-1}c_{12}e_{5}f_{5} e_{5}c_{12}
  x_{6}c_{13}f_{6}^{-1}
c_{14}e_{6}f_{6} \]\medskip

Now, using commutativity and the fact that $t_{f\left( \alpha
\right) }=ft_{\alpha
}f^{-1},$ for any simple closed curve $\alpha $ in $%
\Sigma_{g}$ and any diffeomorphism $f:\Sigma_{g}\rightarrow
\Sigma_{g},$ where  $t_{ \alpha  } $ and $t_{f\left( \alpha \right)
}$ are Dehn twists about the curves $\alpha  $ and $f\left( \alpha
\right) $ respectively, we can write

\[
c_{1}c_{1}c_{2}c_{1}^{-1} x_{1}c_{3}f_{1}^{-1}
c_{4}f_{1}e_{1}e_{1}c_{4}x_{2}c_{5}f_{2}^{-1} c_{6}
  x_{3}c_{7}f_{3}^{-1}c_{8}f_{3} e_{3}e_{3}c_{8}
x_{4}c_{9}f_{4}^{-1} c_{10}
  x_{5}c_{11}f_{5}^{-1}c_{12}f_{5}e_{5} e_{5}c_{12}
  x_{6}c_{13}f_{6}^{-1}
c_{14} \]

as \\
\begin{equation} \label{beginning-part-even} c_{1}d
x_{1}c_{3}r_{1}e_{1}e_{1}c_{4}x_{2}c_{5}f_{2}^{-1} c_{6}
  x_{3}c_{7}r_{3} e_{3}e_{3}c_{8}
x_{4}c_{9}f_{4}^{-1} c_{10}
  x_{5}c_{11}r_{5}e_{5} e_{5}c_{12}
  x_{6}c_{13}f_{6}^{-1}
c_{14}x_{7}c_{15},
\end{equation}

where $d=c_{1}c_{2}c_{1}^{-1}$ and $ r_{i}=f_{i}^{-1}
c_{2i+2}f_{i},i=1,3,5.$\\

The last portion of the word in \eqref{juxtaposed word-even} will be
simplified using the same procedure:

\[\underline{c_{2g-6}e_{g-4}c_{2g-6}}f_{g-4}
x_{g-3}c_{2g-5}e_{g-4}^{-1}f_{g-4}^{-1}e_{g-3}^{-1}f_{g-3}^{-1}
\underline{c_{2g-4}e_{g-3}c_{2g-4}}f_{g-3}\underline{c_{2g-4}e_{g-3}c_{2g-4}}f_{g-3}
    \]\[
 x_{g-2}c_{2g-3}e_{g-3}^{-1}f_{g-3}^{-1}e_{g-2}^{-1}f_{g-2}^{-1} \underline{c_{2g-2}e_{g-2}c_{2g-2}}f_{g-2}
x_{g-1}c_{2g-1}e_{g-2}^{-1}f_{g-2}^{-1}c_{2g+1}^{-1}\underline{c_{2g+1}^{-1}c_{2g+1}}c_{2g}\left(
c_{2g+1}c_{2g}\right)
  ^{3}\]
\[\downarrow\]
\[e_{g-4}c_{2g-6}\underline{e_{g-4}f_{g-4}}
x_{g-3}c_{2g-5}\underline{e_{g-4}^{-1}f_{g-4}}^{-1}\underline{e_{g-3}}^{-1}f_{g-3}^{-1}
\underline{e_{g-3}}c_{2g-4}e_{g-3}f_{g-3}e_{g-3}c_{2g-4}\underline{e_{g-3}f_{g-3}}
    \]\[
 x_{g-2}c_{2g-3}\underline{e_{g-3}^{-1}f_{g-3}}^{-1}\underline{e_{g-2}}^{-1}f_{g-2}^{-1}\underline{ e_{g-2}}c_{2g-2}\underline{e_{g-2}f_{g-2}}
x_{g-1}c_{2g-1}\underline{e_{g-2}^{-1}f_{g-2}}^{-1}c_{2g+1}^{-1}\underline{c_{2g}c_{2g+1}c_{2g}}\left(
c_{2g+1}c_{2g}\right)
  ^{2}\]
\[\downarrow\]
\[c_{2g-6}
x_{g-3}c_{2g-5}f_{g-3}^{-1}
c_{2g-4}\underline{e_{g-3}f_{g-3}}e_{g-3}c_{2g-4}
 x_{g-2}c_{2g-3}f_{g-2}^{-1}    \]\[
c_{2g-2}
x_{g-1}c_{2g-1}\underline{c_{2g+1}^{-1}c_{2g+1}}c_{2g}c_{2g+1}\left(
c_{2g+1}c_{2g}\right)
  ^{2}\]
\[\downarrow\]
\[c_{2g-6}
x_{g-3}c_{2g-5}f_{g-3}^{-1} c_{2g-4}f_{g-3}e_{g-3}e_{g-3}c_{2g-4}
 x_{g-2}c_{2g-3}f_{g-2}^{-1}    \]\[
c_{2g-2} x_{g-1}c_{2g-1}c_{2g}c_{2g+1}\left( c_{2g+1}c_{2g}\right)
  ^{2}\]
\[\downarrow\]
\begin{equation} \label{last-part-even}
 c_{2g-6} x_{g-3}c_{2g-5}r_{g-3}e_{g-3}e_{g-3}c_{2g-4}
 x_{g-2}c_{2g-3}f_{g-2}^{-1}
c_{2g-2} x_{g-1}c_{2g-1}c_{2g}c_{2g+1}c_{2g+1}c_{2g}c_{2g+1}c_{2g},
\end{equation}

where $r_{g-3}=f_{g-3}^{-1} c_{2g-4}f_{g-3}$.\\

 Combining \eqref{beginning-part-even} and \eqref{last-part-even} we obtain

\[c_{1}d x_{1}c_{3}r_{1}e_{1}e_{1}c_{4}x_{2}c_{5}f_{2}^{-1} c_{6}
  x_{3}c_{7}r_{3} e_{3}e_{3}c_{8}
x_{4}c_{9}f_{4}^{-1} c_{10}
  x_{5}c_{11}r_{5}e_{5} e_{5}c_{12}
  x_{6}c_{13}f_{6}^{-1}
c_{14}x_{7}c_{15} \]
\begin{equation} \label{genus-even-achiral-final}
\vdots
\end{equation}
\[c_{2g-6} x_{g-3}c_{2g-5}r_{g-3}e_{g-3}e_{g-3}c_{2g-4}
 x_{g-2}c_{2g-3}f_{g-2}^{-1}
c_{2g-2} x_{g-1}c_{2g-1}c_{2g}c_{2g+1}c_{2g+1}c_{2g}c_{2g+1}c_{2g}\]
\medskip

It seems to be a little challenging to remove the remaining negative
exponents from this last expression at this point. \\

In a more compact form the word can be written as:
\[c_{1}d
x_{1}c_{3}r_{1}e_{1}e_{1}c_{4}x_{2}c_{5}f_{2}^{-1} W_6 W_8\cdots W_g
c_{2g-2}
x_{g-1}c_{2g-1}c_{2g}c_{2g+1}c_{2g+1}c_{2g}c_{2g+1}c_{2g},\]

where  $W_i=c_{2i-6} x_{i-3}c_{2i-5}r_{i-3}e_{i-3}e_{i-3}c_{2i-4}
 x_{i-2}c_{2i-3}f_{i-2}^{-1},i=6,8,\dots,g$.\\

 \subsection{Genus g-odd} Most of the argument will be similar to the even case; we just need to make some
changes on the indices. \\

The following are the words from each component listed
 with alternating signs
  \begin{eqnarray} \label{words-before-gluing}
\nonumber \left( c_{1}c_{2}\right) ^{2} \\
\nonumber \left( c_{4}e_{1}c_{4}f_{1}\right) ^{-1} \\
\nonumber c_{6}e_{2}c_{6}f_{2} \\
 \vdots  \hspace{.5in}  \\
\nonumber \left(c_{2g-2}e_{g-2}c_{2g-2}f_{g-2}\right)^{-1} \\
 \nonumber \left( c_{2g+1}c_{2g}\right) ^{2}.
 \end{eqnarray}

Now, \eqref{d's-and-dbar's-even} becomes

\begin{eqnarray}\label{d's-and-dbar's-odd}
 \nonumber \left( c_{4}e_{1}c_{4}f_{1}\right) ^{-1} =d_2^{-1}\left( c_{4}e_{1}c_{4}f_{1}\right) ^{2}\overline{d}_2^{-1}\\
  \left( c_{8}e_{3}c_{8}f_{3}\right)^{-1}=d_3^{-1}\left( c_{8}e_{3}c_{8}f_{3}\right) ^{2}\overline{d}_3^{-1} \\
\nonumber \vdots  \hspace{.5in}\\
  \nonumber \left(
  c_{2g-2}e_{g-2}c_{2g-2}f_{g-2}\right)^{-1}=d_{g-2}^{-1}\left(
  c_{2g-2}e_{g-2}c_{2g-2}f_{g-2}\right)^{2}\overline{d}_{g-2}^{-1}.
\end{eqnarray}
\eqref{lanterns-solved-for-delta}  and
\eqref{lanterns-solved-for-delta-inverse} are still the same.
Therefore  \eqref{d's-eliminated} becomes

\begin{eqnarray}\label{d's-and-dbar's-eliminated-odd}
 \nonumber \left( c_{4}e_{1}c_{4}f_{1}\right) ^{-1}=x_{1}c_{3}c_{1}^{-1}c_{1}^{-1}e_{1}^{-1}f_{1}^{-1}
 \left( c_{4}e_{1}c_{4}f_{1}\right) ^{2}x_{2}c_{5}e_{1}^{-1}f_{1}^{-1}e_{2}^{-1}f_{2}^{-1}\\
 \nonumber\left( c_{8}e_{3}c_{8}f_{3}\right)^{-1}= x_{3}c_{7}e_{2}^{-1}f_{2}^{-1}e_{3}^{-1}f_{3}^{-1}\left( c_{8}e_{3}c_{8}f_{3}\right) ^{2}
 x_{4}c_{9}e_{3}^{-1}f_{3}^{-1}e_{4}^{-1}f_{4}^{-1} \\
 \vdots  \hspace{1.5in}
\end{eqnarray}
\begin{eqnarray}
  \nonumber \left( c_{2g-2}e_{g-2}c_{2g-2}f_{g-2}\right)^{-1}=x_{g-2}c_{2g-3}e_{g-3}^{-1}f_{g-3}^{-1}e_{g-2}^{-1}f_{g-2}^{-1}\left( c_{2g-2}e_{g-2}c_{2g-2}f_{g-2}\right)^{2}
 \!x_{g-1}c_{2g-1}e_{g-2}^{-1}f_{g-2}^{-1}c_{2g+1}^{-1}c_{2g+1}^{-1}
\end{eqnarray}
using lantern relations in \eqref{d's-and-dbar's-odd} to replace
${d}_i^{-1}$ and $\overline{d}_i^{-1}$.

Then \eqref{words-put-together-even} becomes
\[ \left( c_{1}c_{2}\right) ^{2} \]
\[ x_{1}c_{3}c_{1}^{-1}c_{1}^{-1}e_{1}^{-1}f_{1}^{-1}
 \left( c_{4}e_{1}c_{4}f_{1}\right)
 ^{2}x_{2}c_{5}e_{1}^{-1}f_{1}^{-1}e_{2}^{-1}f_{2}^{-1}\]
\[ c_{6}e_{2}c_{6}f_{2} \]
 \[ x_{3}c_{7}e_{2}^{-1}f_{2}^{-1}e_{3}^{-1}f_{3}^{-1}
 \left( c_{8}e_{3}c_{8}f_{3}\right) ^{2}x_{4}c_{9}e_{3}^{-1}f_{3}^{-1}e_{4}^{-1}f_{4}^{-1}
 \]
\[ c_{10}e_{4}c_{10}f_{4} \]
 \begin{eqnarray} \label{before-juxtaposing-odd}
\vdots  \hspace{.5in}
\end{eqnarray}
\[ c_{2g-4}e_{g-3}c_{2g-4}f_{g-3} \]
\[ x_{g-2}c_{2g-3}e_{g-3}^{-1}f_{g-3}^{-1}e_{g-2}^{-1}f_{g-2}^{-1}\left(
c_{2g-2}e_{g-2}c_{2g-2}f_{g-2}\right)^{2}
 x_{g-1}c_{2g-1}e_{g-2}^{-1}f_{g-2}^{-1}c_{2g+1}^{-1}c_{2g+1}^{-1}\]
\[ \left( c_{2g+1}c_{2g}\right) ^{2}\]

and juxtaposing them results in

\[
 \left( c_{1}c_{2}\right) ^{2}
x_{1}c_{3}c_{1}^{-1}c_{1}^{-1}e_{1}^{-1}f_{1}^{-1}
 \left( c_{4}e_{1}c_{4}f_{1}\right) ^{2}x_{2}c_{5}e_{1}^{-1}f_{1}^{-1}e_{2}^{-1}f_{2}^{-1}
 c_{6}e_{2}c_{6}f_{2} \]\[
  x_{3}c_{7}e_{2}^{-1}f_{2}^{-1}e_{3}^{-1}f_{3}^{-1}\left( c_{8}e_{3}c_{8}f_{3}\right) ^{2}x_{4}c_{9}e_{3}^{-1}f_{3}^{-1}e_{4}^{-1}f_{4}^{-1}
c_{10}e_{4}c_{10}f_{4} \]\[
  x_{5}c_{11}e_{4}^{-1}f_{4}^{-1}e_{5}^{-1}f_{5}^{-1}\left( c_{8}e_{3}c_{8}f_{3}\right) ^{2}x_{6}c_{13}e_{5}^{-1}f_{5}^{-1}e_{6}^{-1}f_{6}^{-1}
c_{14}e_{6}c_{14}f_{6} \]
\begin{eqnarray}\label{juxtaposed word-odd}
\vdots  \hspace{.5in}
\end{eqnarray}
\[c_{2g-4}e_{g-3}c_{2g-4}f_{g-3}
x_{g-2}c_{2g-3}e_{g-3}^{-1}f_{g-3}^{-1}e_{g-2}^{-1}f_{g-2}^{-1}
c_{2g-2}e_{g-2}c_{2g-2}f_{g-2}
    \]
     \begin{eqnarray*}
c_{2g-2}e_{g-2}c_{2g-2}f_{g-2}
x_{g-1}c_{2g-1}e_{g-2}^{-1}f_{g-2}^{-1}c_{2g+1}^{-1}c_{2g+1}^{-1}
c_{2g+1}c_{2g}c_{2g+1}c_{2g}
\end{eqnarray*}
 Eliminating the negative exponents
using braid and commutativity relations in the first half of
\eqref{juxtaposed word-odd} results in
 \begin{eqnarray} \label{juxtaposed word-odd-simplified}\nonumber c_{1}d
x_{1}c_{3}r_{1}e_{1}e_{1}c_{4}x_{2}c_{5}f_{2}^{-1} c_{6}
  x_{3}c_{7}r_{3} e_{3}e_{3}c_{8}
x_{4}c_{9}f_{4}^{-1} c_{10}
  x_{5}c_{11}r_{5}e_{5} e_{5}c_{12}
  x_{6}c_{13}f_{6}^{-1}
c_{14}x_{7}c_{15}, \\
\end{eqnarray}

\noindent where $d=c_{1}c_{2}c_{1}^{-1}$ and $ r_{i}=f_{i}^{-1}
c_{2i+2}f_{i},i=1,3,5,$ just like in the even case ( See \eqref{juxtaposed word-even} - and \eqref{beginning-part-even} ). \\
\indent The last part is also simplified using braid relation first
\[\underline{c_{2g-4}e_{g-3}c_{2g-4}}f_{g-3}
x_{g-2}c_{2g-3}e_{g-3}^{-1}f_{g-3}^{-1}e_{g-2}^{-1}f_{g-2}^{-1}
\underline{c_{2g-2}e_{g-2}c_{2g-2}}f_{g-2}
    \]\[
\underline{c_{2g-2}e_{g-2}c_{2g-2}}f_{g-2}
x_{g-1}c_{2g-1}e_{g-2}^{-1}f_{g-2}^{-1}c_{2g+1}^{-1}c_{2g+1}^{-1}
c_{2g+1}\underline{c_{2g}c_{2g+1}c_{2g}}\]
\[\downarrow\]
\[e_{g-3}c_{2g-4}e_{g-3}f_{g-3}
x_{g-2}c_{2g-3}e_{g-3}^{-1}f_{g-3}^{-1}e_{g-2}^{-1}f_{g-2}^{-1}
e_{g-2}c_{2g-2}e_{g-2}f_{g-2}
    \]\[
e_{g-2}c_{2g-2}e_{g-2}f_{g-2}
x_{g-1}c_{2g-1}e_{g-2}^{-1}f_{g-2}^{-1}c_{2g+1}^{-1}c_{2g+1}^{-1}
c_{2g+1}c_{2g+1}c_{2g}c_{2g+1}\]\\
and commutativity relation along with cancelation next
\[e_{g-3}c_{2g-4}\underline{e_{g-3}f_{g-3}}
x_{g-2}c_{2g-3}\underline{e_{g-3}^{-1}f_{g-3}^{-1}}\
 \underline{e_{g-2}^{-1}}f_{g-2}^{-1} \underline{e_{g-2}}c_{2g-2}e_{g-2}f_{g-2}
    \]\[
e_{g-2}c_{2g-2}\underline{e_{g-2}f_{g-2}}
x_{g-1}c_{2g-1}\underline{e_{g-2}^{-1}f_{g-2}^{-1}}\
 \underline{c_{2g+1}^{-1}c_{2g+1}^{-1}
c_{2g+1}c_{2g+1}}c_{2g}c_{2g+1}\]
\[\downarrow\]
\[e_{g-3}c_{2g-4}
x_{g-2}c_{2g-3}f_{g-2}^{-1}c_{2g-2}e_{g-2}f_{g-2} e_{g-2}c_{2g-2}
x_{g-1}c_{2g-1}
 c_{2g}c_{2g+1}\]
Finally, defining $ r_{g-2}=f_{g-2}^{-1} c_{2g-2}f_{g-2}$ and using
commutativity one more time we obtain:
\[e_{g-3}c_{2g-4}
x_{g-2}c_{2g-3}f_{g-2}^{-1}c_{2g-2}\underline{e_{g-2}f_{g-2}}
e_{g-2}c_{2g-2} x_{g-1}c_{2g-1}
 c_{2g}c_{2g+1}\]
\[\downarrow\]
\[e_{g-3}c_{2g-4}
x_{g-2}c_{2g-3}r_{g-2}f_{g-2}e_{g-2}c_{2g-2} x_{g-1}c_{2g-1}
 c_{2g}c_{2g+1}\]
Putting the two ends together, the word now has the form
\[\nonumber c_{1}d x_{1}c_{3}r_{1}e_{1}e_{1}c_{4}x_{2}c_{5}f_{2}^{-1}
c_{6}
  x_{3}c_{7}r_{3} e_{3}e_{3}c_{8}
x_{4}c_{9}f_{4}^{-1} c_{10}
  x_{5}c_{11}r_{5}e_{5} e_{5}c_{12}
  x_{6}c_{13}f_{6}^{-1}
c_{14}x_{7}c_{15}, \]
\begin{equation}\label{genus-odd-achiral-final}
\vdots
\end{equation}
\[c_{2g-4}
x_{g-2}c_{2g-3}r_{g-2}f_{g-2}e_{g-2}c_{2g-2} x_{g-1}c_{2g-1}
 c_{2g}c_{2g+1}\]
In a more compact form we have
\[c_{1}d x_{1}c_{3}r_{1}e_{1}e_{1}c_{4}x_{2}c_{5}f_{2}^{-1}
 W_6W_8\cdots W_{g-1}  W_g,\]

where $W_i=c_{2i-6} x_{i-3}c_{2i-5}r_{i-3}e_{ i-3}e_{ i-3}c_{2i-4}
 x_{i-2}c_{2i-3}f_{i-2}^{-1}$ and \\
 $W_g=c_{2g-6} x_{g-3}c_{2g-5}r_{g-3}e_{ g-3}e_{ g-3}c_{2g-4}
x_{g-2}c_{2g-3}r_{g-2}f_{g-2}e_{ g-2}c_{2g-2} x_{g-1}c_{2g-1}
 c_{2g}c_{2g+1}$.\\

 The next
section deals with  lower genus.

\section{Low genus} \label{low genus}

For genus 2 and 3 there isn't  much difficulty with eliminating the
terms with negative exponents. For genus $4, 5,$ and $6$ however, we
use   additional lantern relations to eliminate them.

\subsection{genus 2} \label{genus-2-section}

We glue two tori with one one boundary component together and
juxtapose the words $\left( c_{1}c_{2}\right) ^{2}$ and $\left(
c_{5}c_{4}\right) ^{-2}$ on the resulting closed surface.\\

\begin{figure}[htbp]
     \centering  \leavevmode
     \psfig{file=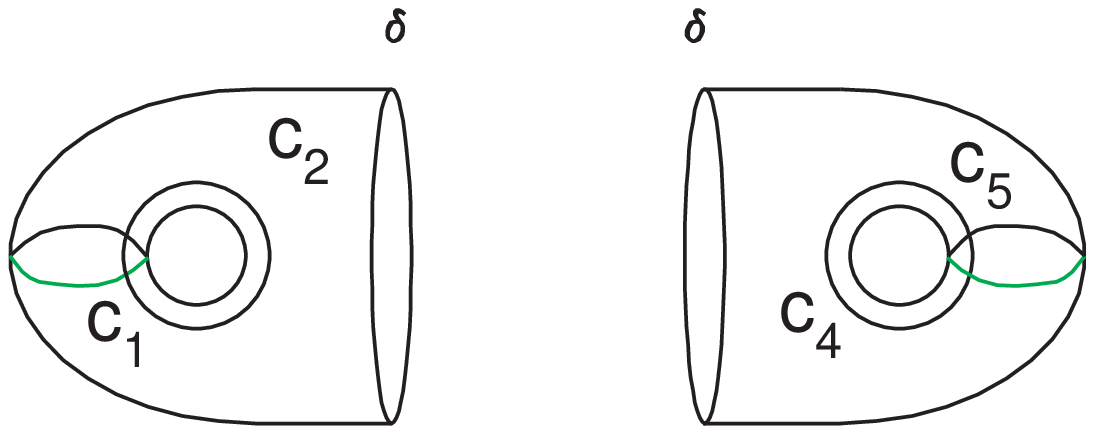, width=5.0in,clip=}
 \end{figure}

\noindent Next, we use the first relation in \eqref{basic relations}
to replace \ $\left( c_{5}c_{4}\right) ^{-2}$  \ by \  $\delta
^{-1}\left( c_{5}c_{4}\right) ^{4}$.\\ Now using the lantern
relation
\[\delta xc_{3}=c_{1}^{2}c_{5}^{2}\] we substitute $\delta ^{-1}=xc_{3}c_{1}^{-2}c_{5}^{-2}$
and obtain
 \[\left( c_{1}c_{2}\right)
^{2}\left( c_{5}c_{4}\right) ^{-2}=\left( c_{1}c_{2}\right)
^{2}\delta ^{-1}\left( c_{5}c_{4}\right) ^{4}=\left(
c_{1}c_{2}\right) ^{2}xc_{3}c_{1}^{-2}c_{5}^{-2}\left(
c_{5}c_{4}\right) ^{4}.\] \vspace{-.2in}
\begin{figure}[htbp]
     \centering  \leavevmode
     \psfig{file=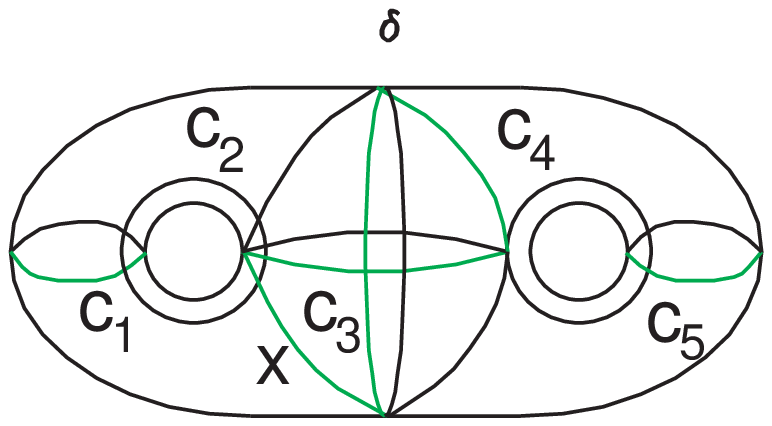, width=6.0in,clip=}
 \end{figure}
\vspace{-.2in}
 Expanding the expression and using
commutativity gives
 \[c_{1}c_{2}c_{1}c_{2}c_{1}^{-2}
xc_{3}c_{5}^{-2}c_{5}c_{4}c_{5}c_{4}c_{5}c_{4}c_{5}c_{4}=c_{1}\underline{c_{2}c_{1}c_{2}}c_{1}^{-2}
xc_{3}c_{5}^{-1}\underline{c_{4}c_{5}c_{4}}c_{5}c_{4}c_{5}c_{4}.\]

Using braid relation on the underlined terms and doing the obvious
cancelations afterward we arrive at
\begin{eqnarray*}c_{1}c_{1}c_{2}c_{1}c_{1}^{-2}
xc_{3}c_{5}^{-1}c_{5}c_{4}c_{5}c_{5}c_{4}c_{5}c_{4}=
c_{1}c_{1}c_{2}c_{1}^{-1} xc_{3}c_{4}c_{5}c_{5}c_{4}c_{5}c_{4}.
\end{eqnarray*}
Even though $c_{1}c_{2}c_{1}^{-1}$ represents a positive twist, we
can do away with this conjugation with little effort: Just bring the
left most $c_1$ in the third power of the word to the right end and
see the cancelations that occur between the underlined terms as you
go from right to left.
\begin{eqnarray} \label{rotation}
\nonumber && \left(c_{1}c_{1}c_{2}c_{1}^{-1}
xc_{3}c_{4}c_{5}c_{5}c_{4}c_{5}c_{4}
\right)^3\\
 \nonumber \ \ \ &=& \!\!\!\!\underline{c_{1}}c_{1}c_{2}c_{1}^{-1}
xc_{3}c_{4}c_{5}c_{5}c_{4}c_{5}c_{4}  c_{1}c_{1}c_{2}c_{1}^{-1}
xc_{3}c_{4}c_{5}c_{5}c_{4}c_{5}c_{4}
c_{1}c_{1}c_{2}c_{1}^{-1} xc_{3}c_{4}c_{5}c_{5}c_{4}c_{5}c_{4} \\
\nonumber \ \ \ &\equiv& \!\!\!\!c_{1}c_{2}\underline{c_{1}^{-1}}
xc_{3}c_{4}c_{5}c_{5}c_{4}c_{5}c_{4}
\underline{c_{1}}c_{1}c_{2}\underline{c_{1}^{-1}}
xc_{3}c_{4}c_{5}c_{5}c_{4}c_{5}c_{4}
\underline{c_{1}}c_{1}c_{2}\underline{c_{1}^{-1}} xc_{3}c_{4}c_{5}c_{5}c_{4}c_{5}c_{4} \underline{c_{1}}\\
\nonumber \ \ \ &=& \!\!\!\!c_{1}c_{2}
xc_{3}c_{4}c_{5}c_{5}c_{4}c_{5}c_{4}\ c_{1}c_{2}
xc_{3}c_{4}c_{5}c_{5}c_{4}c_{5}c_{4}\ c_{1}c_{2}
xc_{3}c_{4}c_{5}c_{5}c_{4}c_{5}c_{4}\\  \ \ \ &=&
\!\!\!\!\left(c_{1}c_{2}
xc_{3}c_{4}c_{5}c_{5}c_{4}c_{5}c_{4}\right)^3=1
\end{eqnarray}

\subsubsection {An Alternate Expression} 
\begin{figure}[htbp]
     \centering  \leavevmode
     \psfig{file=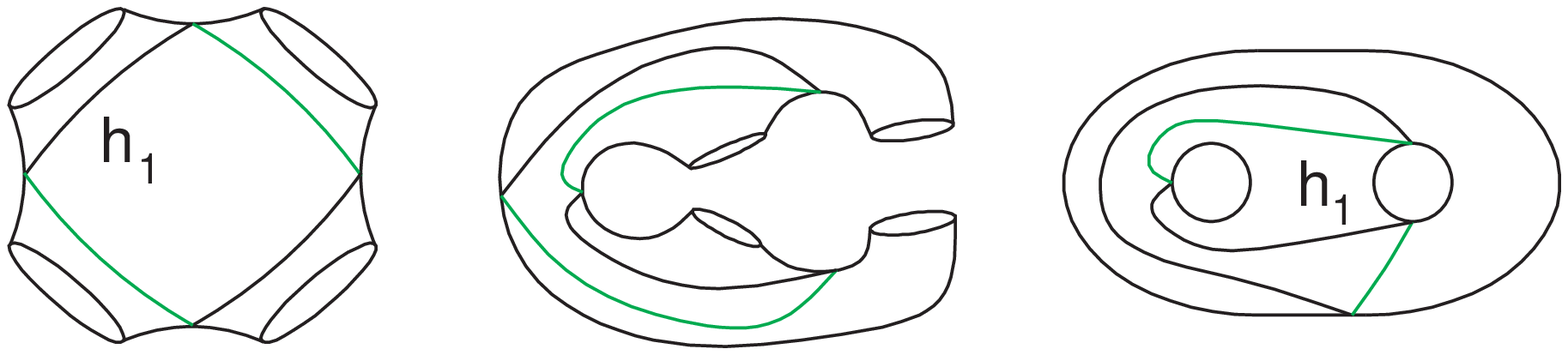, width=5.50in,clip=}
     \caption{ } \label{genus-2-alternate-lantern-figure}
 \end{figure}

We can obtain  an alternate expression out of \eqref{rotation}
by inserting into it  lantern relations as follows:
\begin{eqnarray} \nonumber &&
\left(c_{1}c_{2} xc_{3}c_{4}c_{5}c_{5}c_{4}c_{5}c_{4}\right)^3 \\
\nonumber &\equiv &
\left(c_{1}c_{2} xc_{3}c_{4}c_{3}^{-2}c_{3}^{2}c_{5}^{2}c_{4}c_{5}c_{4}\right)^3 \\
\nonumber &\equiv & \left(c_{1}c_{2} xc_{3}c_{4}c_{3}^{-2}k_1h_1c_{1} c_{4}c_{5}c_{4}\right)^3\\
 &\equiv & \left( c_{2}
xc_{3}c_{4}c_{3}^{-2}k_1h_1c_{1}^{2}c_{5}^{2}c_{5}^{-2}c_{4}c_{5}c_{4}\right)^3 \label{genus-2-lantern-inserted}\\
\nonumber &\equiv & \left( c_{2}
xc_{3}c_{4}c_{3}^{-2}k_1h_1c_{3}\delta
xc_{5}^{-2} c_{4}c_{5}c_{4}\right)^3\\
\nonumber &\equiv & \left( c_{2} xc_{3}c_{4}c_{3}^{-1}k_1 h_1 \delta
xc_{5}^{-2}c_{5} c_{4}c_{5}\right)^3 \\ \nonumber  &\equiv & \left(
c_{2} x\left(c_{3}c_{4}c_{3}^{-1}\right) k_1 h_1 \delta x\left(c_{5}^{-1}c_{4}c_{5}\right)\right)^3 \\
\nonumber  &\equiv & \left( c_{2} xt_2 k_1 h_1 \delta
xs_2\right)^3=1
\end{eqnarray}
where $t_2=c_{3}c_{4}c_{3}^{-1}$ and $s_2=c_{5}^{-1}c_{4}c_{5}$. The
cycles that are used in the first lantern relation
\begin{eqnarray*}
c_1k_1h_1&=&c_3^2c_5^2
\end{eqnarray*} used in  \eqref{genus-2-lantern-inserted} are shown in Figure
\ref{genus-2-alternate-lantern-figure}.\\

\begin{figure}[htbp]
     \centering  \leavevmode
     \psfig{file=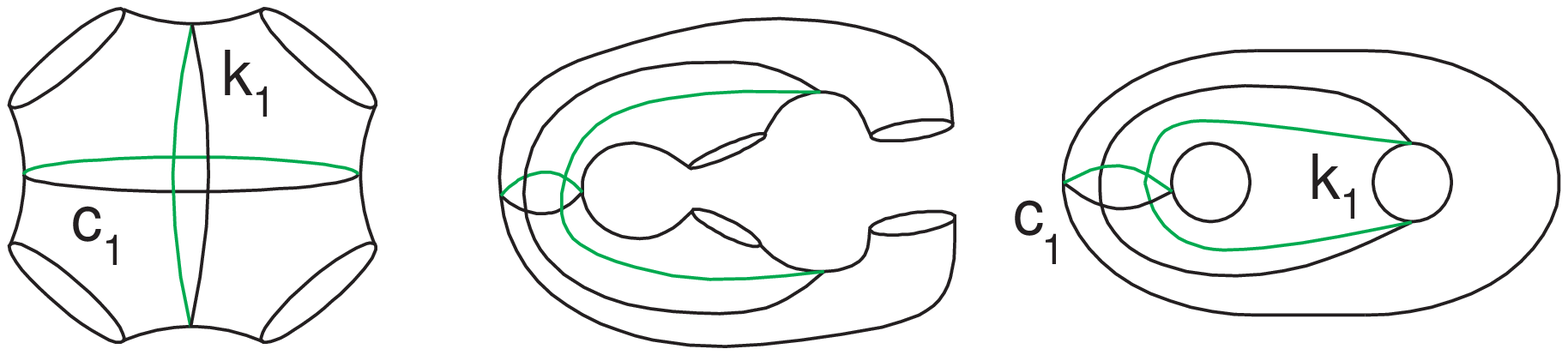, width=5.50in,clip=}
 \end{figure}
 \vspace{-.8in}

\subsection{genus 3} \label{genus-3-section} To get the word for genus 3 we will use \eqref{before-juxtaposing-odd}
with $g=3$:
\begin{eqnarray}
\nonumber \left( c_{1}c_{2}\right) ^{2} \\
\nonumber x_{1}c_{3}c_{1}^{-1}c_{1}^{-1}e_{1}^{-1}f_{1}^{-1}
 \left( c_{4}e_{1}c_{4}f_{1}\right) ^{2}x_{2}c_{5}e_{1}^{-1}f_{1}^{-1}c_{7}^{-1}c_{7}^{-1}\\
\nonumber \left( c_{7}c_{6}\right) ^{2},
 \end{eqnarray}
i.e.
\begin{eqnarray}
\nonumber  c_{1}c_{2} c_{1}c_{2} \
x_{1}c_{3}c_{1}^{-1}c_{1}^{-1}e_{1}^{-1}f_{1}^{-1}
  c_{4}e_{1}c_{4}f_{1} c_{4}e_{1}c_{4}\underline{f_{1}}
 x_{2}c_{5}e_{1}^{-1}\underline{f_{1}^{-1}}c_{7}^{-1}\underline{c_{7}^{-1}}\ \underline{ c_{7}}c_{6}
 c_{7}c_{6}.
 \end{eqnarray}
 After this initial cancelation we use braid relation on the
 underlined terms
\begin{eqnarray}
\nonumber  c_{1}\underline{c_{2} c_{1}c_{2}} \
x_{1}c_{3}c_{1}^{-1}c_{1}^{-1}e_{1}^{-1}f_{1}^{-1}
  \underline{c_{4}e_{1}c_{4}}f_{1} \underline{c_{4}e_{1}c_{4}}
 x_{2}c_{5}e_{1}^{-1}c_{7}^{-1}\underline{c_{6} c_{7}c_{6}},
 \end{eqnarray}
and get
\begin{eqnarray}
\nonumber  c_{1}c_{1}c_{2}\underline{c_{1}} \
x_{1}c_{3}\underline{c_{1}^{-1}}c_{1}^{-1}\underline{e_{1}^{-1}}f_{1}^{-1}
 \underline{ e_{1}}c_{4}e_{1}f_{1} e_{1}c_{4}\underline{e_{1}}
 x_{2}c_{5}\underline{e_{1}^{-1}}\underline{c_{7}^{-1}
 c_{7}}c_{6}c_{7}.
 \end{eqnarray}
Cancelation of the underlined terms gives
\begin{eqnarray}\label{genus-3-word-with-neg-exp}
  c_{1}c_{1}c_{2} \ x_{1}c_{3}c_{1}^{-1}f_{1}^{-1}
 c_{4}e_{1}f_{1} e_{1}c_{4}
 x_{2}c_{5}c_{6}c_{7}.
 \end{eqnarray}
 Now, rearranging the terms using commutativity and letting $r=f_{1}^{-1}
 c_{4}f_{1}$ we obtain
\begin{eqnarray}
\nonumber  c_{1}c_{1}c_{2}c_{1}^{-1} x_{1}c_{3}re_{1} e_{1}c_{4}
 x_{2}c_{5}c_{6}c_{7}.
 \end{eqnarray}
Using the same kind of rotation as in \eqref{rotation} will allow us
to eliminate the conjugation $c_{1}c_{2}c_{1}^{-1}$ and we will get
\begin{eqnarray} \label{genus-3-word-final}
 \left(c_{1}c_{2} x_{1}c_{3}rc_{8} c_{8}c_{4}
 x_{2}c_{5}c_{6}c_{7}\right)^3=1
 \end{eqnarray}
  in the end, using the identification $e_{1}=c_{8}$ as shown in Figure
 \ref{figure-lanterns}.

\subsubsection {An Alternate Expression} An alternate expression is obtained when $f_1^{-1}$ is eliminated
from \eqref{genus-3-word-with-neg-exp} using the lantern relation
\[f_{1}tv=c_{1}c_{3}c_{5}c_{7}.\]
\begin{figure}[htbp]
     \centering  \leavevmode
     \psfig{file=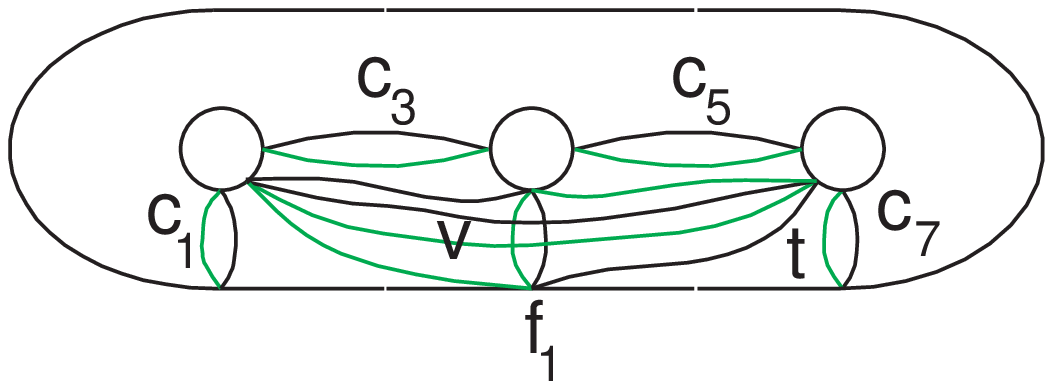,width=5.250in,clip=}
     \caption{} \label{figure-alternate-lanterns-genus-3}
 \end{figure}
Substituting\ \
$f_{1}^{-1}=tvc_{1}^{-1}c_{3}^{-1}c_{5}^{-1}c_{7}^{-1}$ into
\eqref{genus-3-word-with-neg-exp} we get
\[c_{1}c_{1}c_{2} x_{1}\underline{c_{3}}c_{1}^{-1}tvc_{1}^{-1}\underline{c_{3}^{-1}}c_{5}^{-1}c_{7}^{-1}
 c_{4}e_{1}f_{1} e_{1}c_{4}
 x_{2}c_{5}c_{6}c_{7}.\]
We can cancel the underlined terms and rewrite the rest of the word
as
\[c_{1}c_{1}c_{2}c_{1}^{-1}c_{1}^{-1} x_{1}tvc_{5}^{-1}
 c_{4}e_{1}f_{1} e_{1}c_{4}
 x_{2}c_{5}c_{7}^{-1}c_{6}c_{7}\]
using commutativity.  Now, because $c_{5}\left( \bar{x}_{2}\right)
=x_{2}$ ( i.e., Dehn twist about $c_{5}$ maps $\bar{x}_{2}$ to
$x_{2}$ ) we have
\[c_{5}\bar{x}_{2}c_{5}^{-1}=x_{2}, \  \text{i.e.,} \
c_{5}\bar{x}_{2}=x_{2}c_{5}.\] Substituting $c_{5}\bar{x}_{2}$ in
place of $x_{2}c_{5}$ and  inserting a $c_{5}c_{5}^{-1}$  using
commutativity results in
\begin{equation*} \label{genus-3-alternate-word}
c_{1}c_{1}c_{2}c_{1}^{-1}c_{1}^{-1} x_{1}tvc_{5}^{-1}
 c_{4}c_{5}e_{1}f_{1} e_{1}c_{5}^{-1}c_{4}
c_{5}\bar{x}_{2}c_{7}^{-1}c_{6}c_{7}.
\end{equation*}
\begin{figure}[htbp]
     \centering  \leavevmode
     \psfig{file=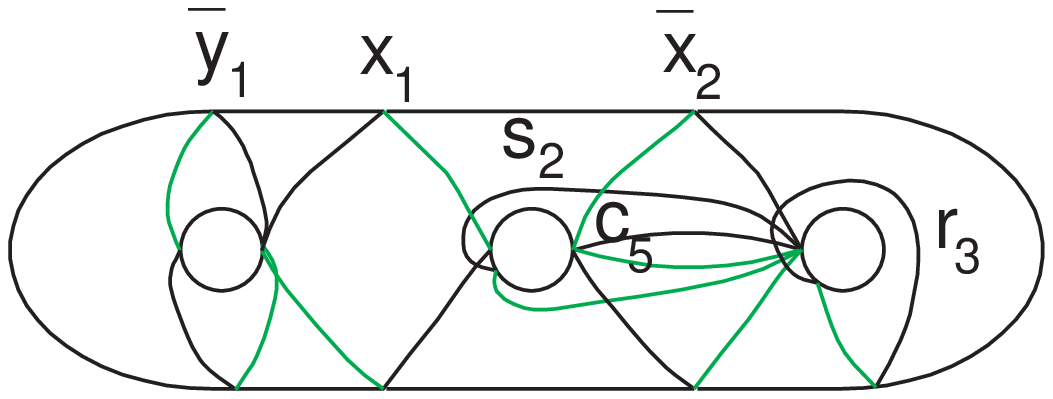,width=5.250in,clip=}
     \caption{} \label{figure-alternate-cycles-genus-3}
 \end{figure}
 Now, all we have to do is rename the conjugations. If we let
 \[\bar{y}_1=c_{1}c_{1}c_{2}c_{1}^{-1}c_{1}^{-1},s_2=c_{5}^{-1}
 c_{4}c_{5},\ \  \text{and} \ \ r_3=c_{7}^{-1}c_{6}c_{7}\]
 then the final form of the word becomes
 \begin{equation} \label{genus-3-alternate-word-final}
 \left(\bar{y}_1 x_{1}tvs_2c_{8}f_{1}
 c_{8}s_2\bar{x}_{2}r_3\right)^3=1,
 \end{equation}
using the identification $e_1=c_8$ in Figure \ref{figure-lanterns}.

\subsubsection {Another Alternate Expression} One other  expression is obtained using
the relation $$\left( c_{1}c_{2}c_{3}\right) ^{4}=e_{1}f_{1}$$ in
order to substitute $\left( c_{1}c_{2}c_{3}\right) ^{4}$ in place of
$c_{8}f_{1}=e_{1}f_{1}$. The alternate expression is
\begin{equation} \label{genus-3-second-alternate-word}
 \left(\bar{y}_1 x_{1}tvs_2\left( c_{1}c_{2}c_{3}\right) ^{4}
 c_{8}s_2\bar{x}_{2}r_3\right)^3=1,
 \end{equation}

 \subsection{genus 4}\label{genus-4-section} Using \eqref{words-put-together-even} with
 $g=4$ we have

\[ \left( c_{1}c_{2}\right) ^{2} \]
\[ x_{1}c_{3}c_{1}^{-1}c_{1}^{-1}e_{1}^{-1}f_{1}^{-1}
 \left( c_{4}e_{1}c_{4}f_{1}\right)
 ^{2}x_{2}c_{5}e_{1}^{-1}f_{1}^{-1}e_{2}^{-1}f_{2}^{-1}\]
 \[ c_{6}e_{2}c_{6}f_{2} \]
 \[ x_{3}c_{7}e_{2}^{-1}f_{2}^{-1}c_{9}^{-1}c_{9}^{-1}\left(
c_{9}c_{8}\right)
  ^{4},\]
which will be the same as the top line in
\eqref{beginning-of-the-juxtaposed-word} for the most part after
juxtaposing:
\[c_{1}c_{1}c_{2}\underline{c_{1}}
x_{1}c_{3}\underline{c_{1}}^{-1}c_{1}^{-1}\underline{e_{1}}^{-1}f_{1}^{-1}
 \underline{e_{1}}c_{4}e_{1}f_{1}e_{1}c_{4}\underline{e_{1}f_{1}}x_{2}c_{5}\underline{e_{1}^{-1}f_{1}}^{-1}\underline{e_{2}}^{-1}f_{2}^{-1}
\underline{e_{2}}c_{6} \underline{e_{2}f_{2}}x_{3}c_{7}\ \
\underline{e_{2}^{-1}f_{2}^{-1}}c_{9}^{-1}\underline{c_{9}^{-1}}\left(
\underline{c_{9}}c_{8}\right)^4.\] Cancelation of the underlined
terms gives
\[c_{1}c_{1}c_{2}
x_{1}c_{3}c_{1}^{-1}f_{1}^{-1}
c_{4}e_{1}f_{1}e_{1}c_{4}x_{2}c_{5}f_{2}^{-1} c_{6} x_{3}c_{7}\ \
c_{9}^{-1}\underline{c_{8}c_{9}c_{8}} c_{9}c_{8}c_{9}c_{8}.\]
Rearranging some commuting terms along with  braid relation on the
underlined triple we get
\[c_{1}c_{1}c_{2}c_{1}^{-1}
x_{1}c_{3}f_{1}^{-1} c_{4}f_{1}e_{1}e_{1}c_{4}x_{2}c_{5}f_{2}^{-1}
c_{6} x_{3}c_{7}\underline{c_{9}^{-1}c_{9}}c_{8}c_{9}
c_{9}c_{8}c_{9}c_{8},\] and by setting $
d=c_{1}c_{2}c_{1}^{-1},r_1=f_{1}^{-1} c_{4}f_{1}$ and canceling the
underlined pair we arrive at
\begin{equation}\label{genus-4-word-with-negative}
c_{1}d x_{1}c_{3}r_1e_{1}e_{1}c_{4}x_{2}c_{5}f_{2}^{-1} c_{6}
x_{3}c_{7}c_{8}c_{9} c_{9}c_{8}c_{9}c_{8}.
\end{equation}

\begin{figure}[htbp]
     \centering  \leavevmode
     \psfig{file=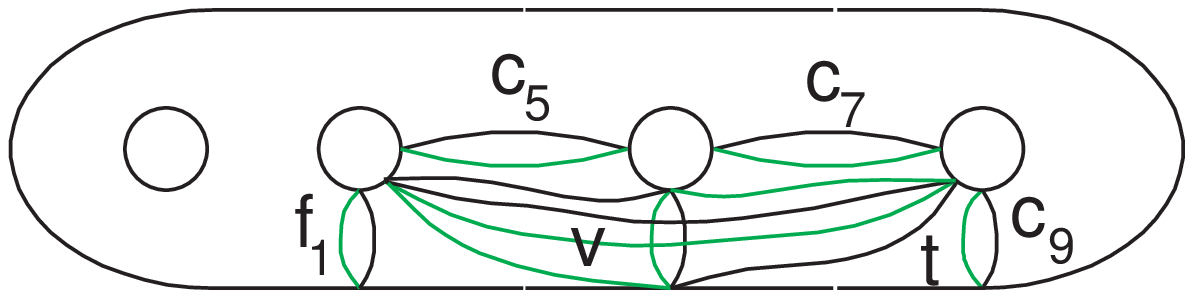,width=5.50in,clip=}
     \caption{} \label{figure-alternate-lanterns-genus-4}
 \end{figure}

We will have to use another lantern relation to eliminate
$f_{2}^{-1}$ and that will be
\[f_2tv=f_1c_5c_7c_9\]
as shown in Figure \ref{figure-alternate-lanterns-genus-4}. We solve
it for $f_2^{-1}$
\[f_2^{-1}=tvf_1^{-1}c_5^{-1}c_7^{-1}c_9^{-1}\]

\noindent and substituting that in
\eqref{genus-4-word-with-negative} using commutativity gives :
\begin{equation*}
c_{1}d
x_{1}c_{3}r_1e_{1}e_{1}c_{4}f_1^{-1}x_{2}\underline{c_{5}}tv\underline{c_5}^{-1}c_7^{-1}
c_{6} x_{3}c_{7}c_9^{-1}c_{8}c_{9}\underline{ c_{9}c_{8}c_{9}}c_{8}.
\end{equation*}
Because $c_{7}\left( \bar{x}_{3}\right) =x_{3}$ ( i.e., Dehn twist
about $c_{7}$ maps $\bar{x}_{3}$ to $x_{3}$ ) we have
\begin{equation}\label{x-conjugate}
c_{7}\bar{x}_{3}c_{7}^{-1}=x_{3}, \  \text{i.e.,} \
c_{7}\bar{x}_{3}=x_{3}c_{7}.
\end{equation}

 Substituting that in and using braid
relation and cancelation on the underlined parts we get
\begin{equation*}
c_{1}d x_{1}c_{3}r_1e_{1}e_{1}c_{4}f_1^{-1}x_{2}tv c_7^{-1} c_{6}
c_{7}\overline{x}_{3}c_9^{-1}\underline{c_{8}c_{9}
c_{8}}c_{9}c_{8}c_{8}.
\end{equation*} Renaming $c_7^{-1} c_{6} c_{7}=s_3$ and using braid
relation on the underlined part again gives
\begin{equation}\label{genus-4-second-from-last}
c_{1}d x_{1}c_{3}\underline{r_1e_{1}e_{1}c_{4}f_1}^{-1}x_{2}tv s_{3}
\overline{x}_{3}\underline{c_9^{-1}c_{9}}c_{8} c_{9}c_{9}c_{8}c_{8}.
\end{equation}
 The following is how we eliminate $f_1^{-1}$ from the underlined
portion:
\begin{eqnarray} \label{from-r-to-y} \nonumber
\underline{r_1}e_{1}e_{1}c_{4}f_1^{-1}= f_{1}^{-1}
c_{4}f_{1}e_{1}e_{1}c_{4}f_1^{-1}=\\  \nonumber f_{1}^{-1}
c_{4}e_{1}e_{1}\underline{f_{1}c_{4}f_1^{-1}}=\\ \nonumber
f_{1}^{-1} c_{4}e_{1}e_{1}c_4^{-1}f_{1}c_{4}=\\ \nonumber f_{1}^{-1}
c_{4}e_{1}\underline{c_4^{-1}f_{1}f_{1}^{-1}
c_{4}}e_{1}c_4^{-1}f_{1}c_{4}=\\  f_{1}^{-1} c_{4}e_{1}c_4^{-1}f_{1}
f_{1}^{-1} c_{4}e_{1}c_4^{-1}f_{1}c_{4}=y_2y_2c_4,
\end{eqnarray}
where $y_2=f_{1}^{-1} c_{4}e_{1}c_4^{-1}f_{1}$. Now
\eqref{genus-4-second-from-last} becomes
\begin{equation*}
c_{1}d x_{1}c_{3}y_2y_2c_4x_{2}tv s_{3} \overline{x}_{3}c_{8}
c_{9}c_{9}c_{8}c_{8}.
\end{equation*}
Using the same rotation operation as in \eqref{rotation}  allows us
to eliminate the conjugation $d=c_{1}c_{2}c_{1}^{-1}$ and we  get
\begin{equation} \label{final-word-genus-4}
\left(c_{1}c_2 x_{1}c_{3}y_2y_2c_4x_{2}tv s_{3}
\overline{x}_{3}c_{8} c_{9}c_{9}c_{8}c_{8}\right)^{3}=1.
\end{equation}

\subsubsection {An Alternate Expression} An alternate expression is
obtained when $f_1^{-1}$ is eliminated from
\eqref{genus-4-second-from-last} using the lantern relation
\[f_{1}t_{1,4}v_{1,4}=c_{1}c_{3}v c_{9}.\]
\begin{figure}[htbp]
     \centering  \leavevmode
     \psfig{file=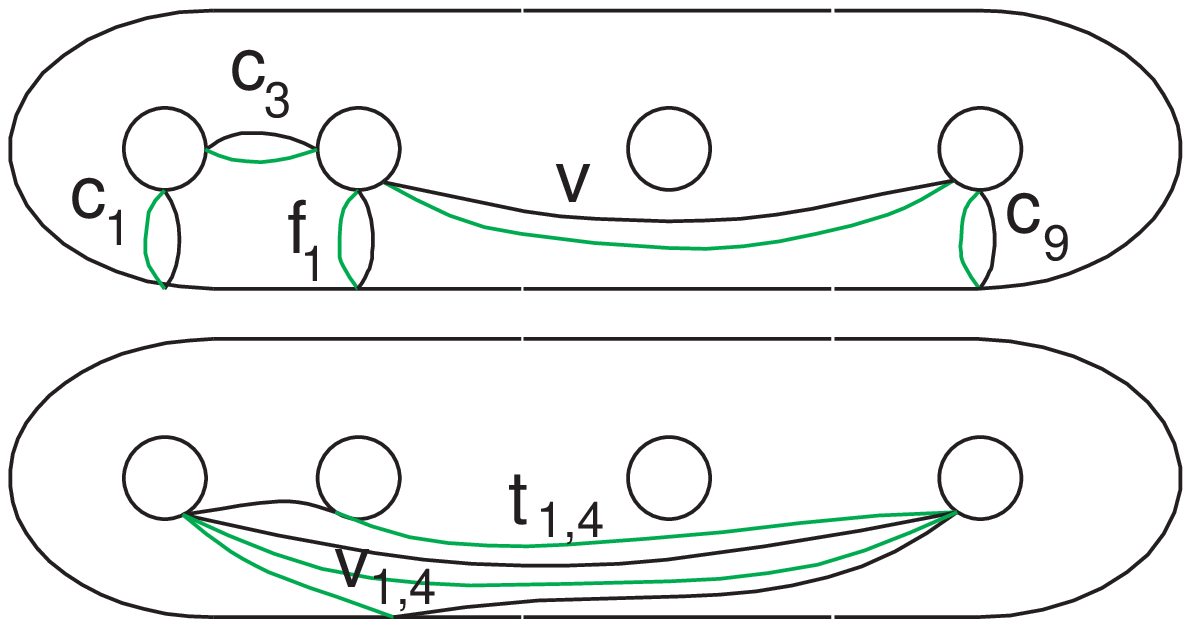,width=5.50in,clip=}
     \caption{}\label{figure-genus-4}
 \end{figure}
Substituting\ \ $f_{1}^{-1}=
t_{1,4}v_{1,4}c_{1}^{-1}c_{3}^{-1}v^{-1}c_{9}^{-1}$ into
\eqref{genus-4-second-from-last} yields
\begin{equation}\label{genus-4-with-second-lantern}
c_{1}d
x_{1}c_{3}r_1e_{1}e_{1}c_{4}t_{1,4}v_{1,4}c_{1}^{-1}c_{3}^{-1}v^{-1}c_{9}^{-1}x_{2}tv
s_{3} \overline{x}_{3}c_{8} c_{9}c_{9}c_{8}c_{8}.
\end{equation}
Using commutativity and inserting identity where necessary we can
rewrite the last expression as
\begin{equation}\label{genus-4-with-second-lantern}
c_{1}c_{1}c_{2}c_{1}^{-1}c_{1}^{-1}
x_{1}c_{3}r_1c_{3}^{-1}e_{1}e_{1}c_{3}c_{4}c_{3}^{-1}t_{1,4}v_{1,4}v^{-1}x_{2}vv^{-1}tv
s_{3} \overline{x}_{3}c_{9}^{-1}c_{8} c_{9}c_{9}c_{8}c_{8},
\end{equation}
remembering $d=c_{1}c_{2}c_{1}^{-1}$.
 Now, all we have to do is rename the conjugations:
 \[\bar{y}_1=c_{1}c_{1}c_{2}c_{1}^{-1}c_{1}^{-1},u_1=c_{3}r_1c_{3}^{-1},\bar{s_2}=c_{3}c_{4}c_{3}^{-1},
 w=v^{-1}x_{2}v,z=v^{-1}tv\ \  \text{and} \ \ r_4=c_{9}^{-1}c_{8}c_{9}.\]
 Then  \eqref{genus-4-with-second-lantern} becomes $\bar{y}_1x_{1}u_1e_{1}e_{1}\bar{s_2}t_{1,4}v_{1,4}wz s_{3} \overline{x}_{3}
r_4c_{9}c_{8}c_{8}$. Therefore the final form of the alternate word
for genus 4 is
\begin{equation}\label{genus-4-with-second-lantern-last}
\left(\bar{y}_1x_{1}u_1e_{1}e_{1}\bar{s_2}t_{1,4}v_{1,4}wz s_{3}
\overline{x}_{3} r_4c_{9}c_{8}c_{8}\right)^3=1.
\end{equation}

\subsubsection {An alternate Construction} \label{An alternate construction}
An alternate gluing operation for genus 4 can be performed as shown
in Figure \ref{figure-genus-4}.

\begin{figure}[htbp]
     \centering  \leavevmode
     \psfig{file=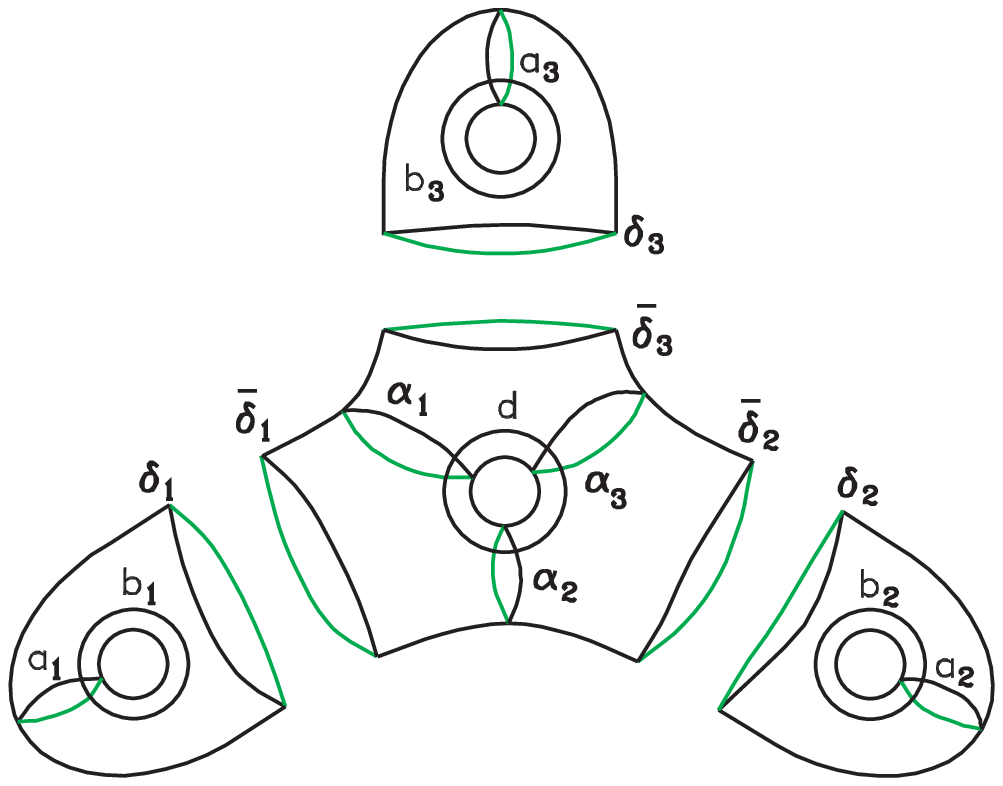,width=6.250in,clip=}
     \caption{}\label{figure-genus-4}
 \end{figure}

We use the words
\begin{eqnarray}\label{genus-4-alternate-word}
\nonumber b _{1}a_1b _{1}a _{1}\\
\nonumber b_{2}a _{2}b_{2}a _{2}\\
b_{3}a_{3}b_{3}a_{3}\\
\nonumber \left( \alpha _{1}\alpha _{2}\alpha _{3} d \right) ^{-1}
\end{eqnarray} on the four
bounded surfaces taking the one in the center with the opposite
orientation. Using the star relation \eqref{basic relations}
\[\left(
\alpha _{1}\alpha _{2}\alpha _{3} d \right)
^{3}=\bar{\delta}_{1}\bar{\delta}_{2}\bar{\delta}_{3},\] we write
\[\left( \alpha _{1}\alpha _{2}\alpha _{3} d\right)
^{-1}=\bar{\delta}_{1}^{-1}\bar{\delta}_{2}^{-1}\bar{\delta}_{3}^{-1}\left(
\alpha _{1}\alpha _{2}\alpha _{3} d \right) ^{2}\] and using the
lantern relations
\begin{eqnarray*}
\delta _{1}x_{1}c_{1}=a_{1}a_{1}\alpha _{1}\alpha _{2}\\
\delta _{2}x_{2}c_{2}=a_{2}a_{2}\alpha _{2}\alpha _{3}\\
\delta _{3}x_{3}c_{3}=a_{3}a_{3}\alpha _{3}\alpha _{1}
\end{eqnarray*}
and the fact that $\delta _{1}=\bar{\delta}_{1},\delta
_{2}=\bar{\delta}_{2},\delta _{3}=\bar{\delta}_{3}$ we write
\begin{figure}[htbp] \label{genus-4}
     \centering  \leavevmode
     \psfig{file=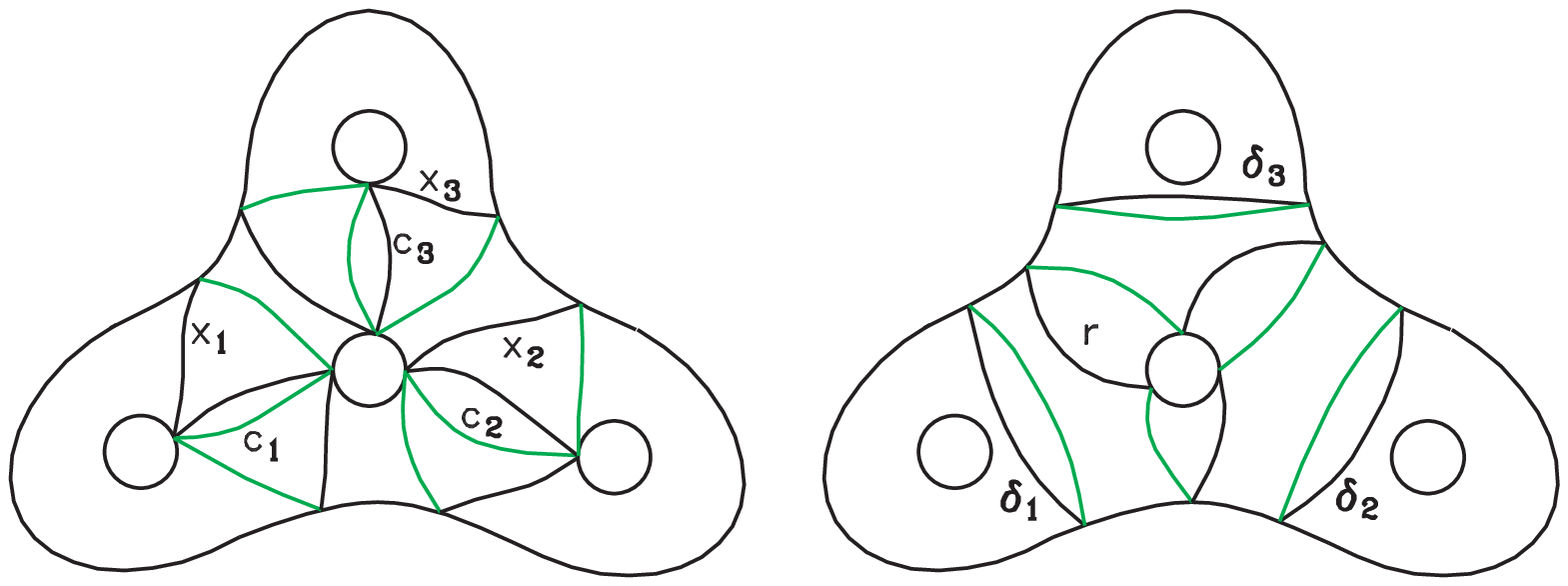,width=6.250in,clip=}
     \caption{}
 \end{figure}
 \begin{eqnarray*}
\delta _{1}^{-1}=x_{1}c_{1}a_{1}^{-1}a_{1}^{-1} \alpha
_{1}^{-1}\alpha _{2}^{-1}\\
\delta _{2}^{-1}=x_{2}c_{2}a_{2}^{-1}a_{2}^{-1} \alpha
_{2}^{-1}\alpha_{3}^{-1}\\
\delta _{3}^{-1}=x_{3}c_{3}a_{3}^{-1}a_{3}^{-1}\alpha
_{3}^{-1}\alpha _{1}^{-1}
\end{eqnarray*}
 Substituting all these in
\eqref{genus-4-alternate-word} and juxtaposing we obtain
\[b _{1}a _{1}b_{1}\underline{a _{1}}b _{2}a _{2}b _{2}\underline{a _{2}}
b_{3}a_{3}b_{3}\underline{a _{3}}x_{1}c_{1}\underline{a _{1}^{-1}}a
_{1}^{-1}\alpha_{1}^{-1}\alpha_{2}^{-1}x_{2}c_{2}\underline{a
_{2}^{-1}}a _{2}^{-1}\alpha_{2}^{-1}\alpha_{3}^{-1}
x_{3}c_{3}\underline{a _{3}^{-1}}a
_{3}^{-1}\alpha_{3}^{-1}\alpha_{1}^{-1} \left(
\alpha_{1}\alpha_{2}\alpha_{3}d\right) ^{2}\] We can cancel the
underlined terms right away using commutativity and rearrange rest
of the word as
\[\underline{b _{1}a _{1}b _{1}}a
_{1}^{-1}\underline{b _{2}a _{2}b _{2}}a _{2}^{-1} \underline{b
_{3}a_{3}b _{3}}a _{3}^{-1}x_{1}c_{1}x_{2}c_{2}
x_{3}c_{3}\alpha_{1}^{-1}\alpha_{2}^{-1}\underline{\alpha_{2}^{-1}}\alpha_{3}^{-1}\underline{\alpha_{3}^{-1}}\
\underline{\alpha_{1}^{-1}}\  \underline{ \alpha_{1}}\
\underline{\alpha_{2}}\
\underline{\alpha_{3}}d\alpha_{1}\alpha_{2}\alpha_{3}d\] using
commutativity again. Now, using braid relation and cancelation on
the underlined portion, the word reduces to
\[a _{1}b _{1}\underline{a _{1}}\ \underline{a
_{1}^{-1}}a _{2}b _{2}\underline{a _{2}}\ \underline{a _{2}^{-1}
}a_{3}b _{3}\underline{a _{3}} \ \underline{a
_{3}^{-1}}x_{1}c_{1}x_{2}c_{2}
x_{3}c_{3}\alpha_{1}^{-1}\alpha_{2}^{-1}\alpha_{3}^{-1}
d\alpha_{1}\alpha_{2}\alpha_{3}d\] Further cancelation and renaming
$r=\left(\alpha_{1}\alpha_{2}\alpha_{3}\right)^{-1}d\alpha_{1}\alpha_{2}\alpha_{3}$
gives the positive relator
\begin{equation} \label{genus-4-rose-word-final}
\left(a _{1}b _{1}a _{2}b _{2}a_{3}b _{3}x_{1}c_{1}x_{2}c_{2}
x_{3}c_{3}rd\right)^3=1
\end{equation}

\subsubsection {Another alternate expression} We can modify \eqref{genus-4-rose-word-final} in order to
insert the lantern relation
\begin{equation*}
 \alpha_{2}t v =a _{1}c_{1}a _{2}c_{2}
\end{equation*}
\begin{figure}[htbp] \label{genus-4}
     \centering  \leavevmode
     \psfig{file=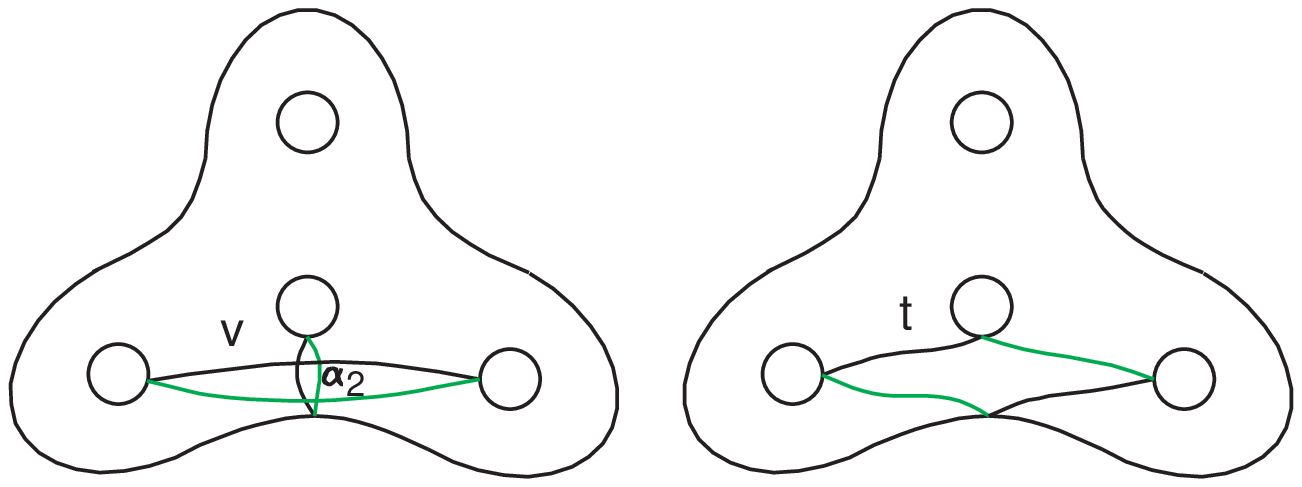,width=6.250in,clip=}
     \caption{}
 \end{figure}
 into it and obtain a new expression.
\begin{eqnarray} \label{genus-4-rose-word-modified}
\nonumber && \left(a _{1}b _{1}a _{1}^{-1}a _{2}b _{2}a
_{2}^{-1}a_{3}b _{3}x_{1}x_{2}a _{1}a _{2}c_{1}c_{2}
x_{3}c_{3}rd\right)^3\\ &=& \left(g_1g_2a_{3}b _{3}x_{1}x_{2}
\alpha_{2}tv x_{3}c_{3}rd\right)^3=1,
\end{eqnarray}
 where $g_1=a _{1}b _{1}a
_{1}^{-1} $ and $g_2=a _{2}b _{2}a _{2}^{-1}$.
 \subsection{genus 5} Using \eqref{before-juxtaposing-odd} with
 $g=5$ we have
 \begin{eqnarray}
\nonumber \left( c_{1}c_{2}\right) ^{2} \\
\nonumber x_{1}c_{3}c_{1}^{-1}c_{1}^{-1}e_{1}^{-1}f_{1}^{-1}
 \left( c_{4}e_{1}c_{4}f_{1}\right) ^{2}x_{2}c_{5}e_{1}^{-1}f_{1}^{-1}e_{2}^{-1}f_{2}^{-1}\\
\nonumber c_{6}e_{2}c_{6}f_{2} \\
 \nonumber x_{3}c_{7}e_{2}^{-1}f_{2}^{-1}e_{3}^{-1}f_{3}^{-1}
 \left( c_{8}e_{3}c_{8}f_{3}\right) ^{2}x_{4}c_{9}e_{3}^{-1}f_{3}^{-1}c_{11}^{-1}c_{11}^{-1} \\
\nonumber c_{11}c_{10}c_{11}c_{10}
 \end{eqnarray}
The top two lines in \eqref{juxtaposed word-odd} up to $c_{10}$
becomes
\begin{eqnarray*} \nonumber c_{1}d
x_{1}c_{3}r_{1}e_{1}e_{1}c_{4}x_{2}c_{5}f_{2}^{-1} c_{6}
  x_{3}c_{7}r_{3} e_{3}e_{3}c_{8}
x_{4}c_{9}c_{11}^{-1} c_{11}^{-1} \\
\end{eqnarray*}
in \eqref{juxtaposed word-odd-simplified}, and this followed by
\[c_{11}c_{10}c_{11}c_{10}=c_{11}c_{11}c_{10}c_{11}\]
gives
\begin{eqnarray*} \nonumber c_{1}d
x_{1}c_{3}r_{1}e_{1}e_{1}c_{4}x_{2}c_{5}f_{2}^{-1} c_{6}
  x_{3}c_{7}r_{3} e_{3}e_{3}c_{8}
x_{4}c_{9}c_{11}^{-1} c_{11}^{-1}c_{11}c_{11}c_{10}c_{11}=
\end{eqnarray*}
\begin{eqnarray} \label{additional-lantern-genus-5}
c_{1}d
x_{1}c_{3}r_{1}e_{1}e_{1}c_{4}x_{2}c_{5}f_{2}^{-1} c_{6}
  x_{3}c_{7}r_{3} e_{3}e_{3}c_{8}
x_{4}c_{9}c_{10}c_{11}
\end{eqnarray}

\begin{figure}[htbp] \label{figure-alternate-lanterns-genus-5}
     \centering  \leavevmode
     \psfig{file=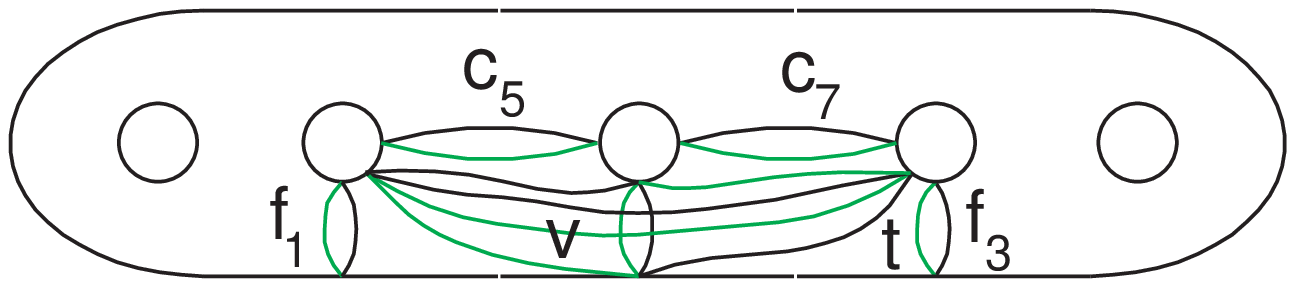,width=5.50in,clip=}
     \caption{}
 \end{figure}

We will use the same additional lantern relation as in genus 4 to
eliminate $f_2^{-1}$: Substitute
\[f_2^{-1}=tvf_1^{-1}c_5^{-1}c_7^{-1}f_3^{-1}\]
in \eqref{additional-lantern-genus-5} using commutativity
\begin{eqnarray*}
c_{1}d
x_{1}c_{3}r_{1}e_{1}e_{1}c_{4}f_1^{-1}x_{2}c_{5}tvc_5^{-1}c_7^{-1}
c_{6}  x_{3}c_{7}f_3^{-1}r_{3} e_{3}e_{3}c_{8}
x_{4}c_{9}c_{10}c_{11}
\end{eqnarray*}

\noindent and eliminate $f_1^{-1},c_5^{-1}$, and $c_7^{-1}$  in the
exact same way as in genus 4 following Figure
\ref{figure-alternate-lanterns-genus-4}. Therefore we can borrow the
portion of \eqref{final-word-genus-4} up to $c_8$ and write
\begin{eqnarray} \label{genus-5-with-one-negative}
c_{1}d x_{1}c_{3}y_2y_2c_4x_{2}tv s_{3}
\overline{x}_{3}\underline{f_3^{-1}r_{3} e_{3}e_{3}c_{8}}
x_{4}c_{9}c_{10}c_{11}.
\end{eqnarray}
The following is how we deal with $f_3^{-1}$:
\begin{eqnarray*}
f_3^{-1}r_{3}
e_{3}e_{3}c_{8}=f_3^{-1}\underline{f_3^{-1}c_8f_3}e_{3}e_{3}c_{8}=\\
f_3^{-1}c_8f_3c_8^{-1}e_{3}e_{3}c_{8}=f_3^{-1}c_8f_3c_8^{-1}e_{3}c_8c_8^{-1}e_{3}c_{8}=r_3\overline{r}_3
\overline{r}_3,
\end{eqnarray*}where
$r_3=f_3^{-1}c_8f_3,\overline{r}_3=c_8^{-1}e_{3}c_8$. Therefore the
final form of the genus 5 word is
\begin{eqnarray}\label{genus-5-word-final}
\left(c_{1}c_2 x_{1}c_{3}y_2y_2c_4x_{2}tv s_{3}
\overline{x}_{3}r_3\overline{r}_3 \overline{r}_3
x_{4}c_{9}c_{10}c_{11}\right)^3=1
\end{eqnarray}
after a rotation similar to \eqref{rotation} applied.

 \subsection{genus 6} Setting $g=6$ in
 \eqref{words-put-together-even} we obtain the  components
\[ \left( c_{1}c_{2}\right) ^{2} \]
\[ x_{1}c_{3}c_{1}^{-1}c_{1}^{-1}e_{1}^{-1}f_{1}^{-1}
 \left( c_{4}e_{1}c_{4}f_{1}\right)
 ^{2}x_{2}c_{5}e_{1}^{-1}f_{1}^{-1}e_{2}^{-1}f_{2}^{-1}\]
 \[ c_{6}e_{2}c_{6}f_{2} \]
 \[ x_{3}c_{7}e_{2}^{-1}f_{2}^{-1}e_{3}^{-1}f_{3}^{-1}\left( c_{8}e_{3}c_{8}f_{3}\right) ^{2}x_{4}c_{9}e_{3}^{-1}f_{3}^{-1}e_{4}^{-1}f_{4}^{-1}
 \]
\begin{equation}c_{10}e_{4}c_{10}f_{4} \label{genus-6-first-stage}
\end{equation}
\[x_{5}c_{11}e_{4}^{-1}f_{4}^{-1}c_{13}^{-1}c_{13}^{-1}\left(
c_{13}c_{12}\right)
  ^{4}\]
of the word  on
 bounded subsurfaces before juxtaposition. After juxtaposing them we arrive at
 \[c_{1}d x_{1}c_{3}r_{1}e_{1}e_{1}c_{4}x_{2}c_{5}f_{2}^{-1} c_{6}
  x_{3}c_{7}r_{3} e_{3}e_{3}c_{8}
x_{4}c_{9}f_{4}^{-1} c_{10}
  x_{5}c_{11}c_{12}c_{13}c_{13}c_{12}c_{13}c_{12}\]
  as in  \eqref{genus-even-achiral-final} with  $g=6$.
If we substitute
\[f_2^{-1}=t_{2,4}v_{2,4}f_1^{-1}c_5^{-1}c_7^{-1}f_3^{-1}\ \ \text{and}\  \ f_4^{-1}=t_{4,6}v_{4,6}f_3^{-1}c_9^{-1}c_{11}^{-1}c_{13}^{-1}\]

\begin{figure}[htbp] \label{figure-alternate-lanterns-f2-f4-genus-6}
     \centering  \leavevmode
     \psfig{file=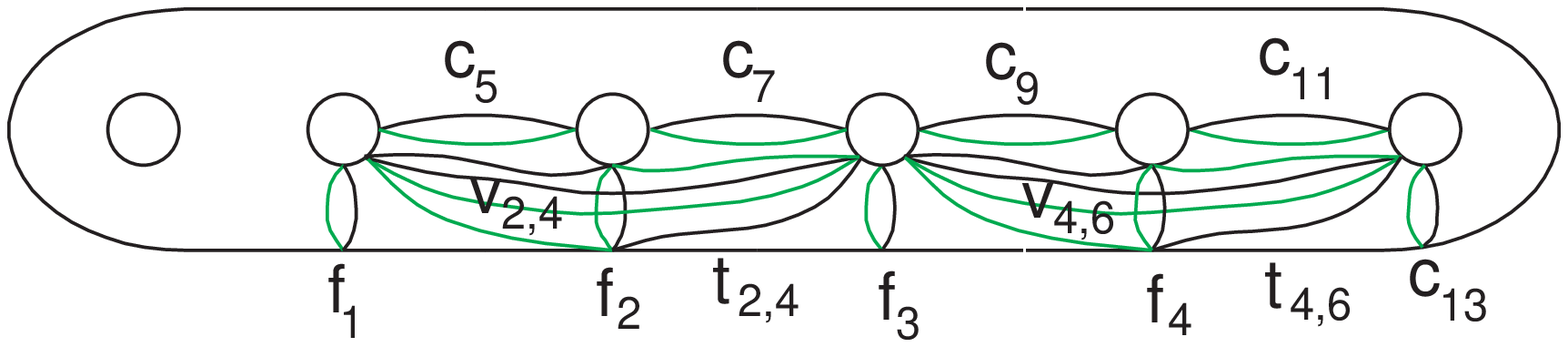,width=5.50in,clip=}
     \caption{}
 \end{figure}
 then we get
\[\hspace{-.5in}c_{1}d x_{1}c_{3}r_{1}e_{1}e_{1}c_{4}x_{2}\underline{c_{5}}t_{2,4}v_{2,4}f_1^{-1}\underline{c_5}^{-1}c_7^{-1}f_3^{-1} c_{6}
  x_{3}c_{7}r_{3} e_{3}e_{3}c_{8}
x_{4}\underline{c_{9}}t_{4,6}v_{4,6}f_3^{-1}\underline{c_9}^{-1}c_{11}^{-1}c_{13}^{-1}c_{10}
  x_{5}c_{11}c_{12}c_{13}c_{13}c_{12}c_{13}c_{12}.\]
  We do the obvious cancelations and use \eqref{x-conjugate} again
  to write
  \[ x_{3}c_{7}=c_{7}\bar{x}_{3} \ \ \text{and likewise} \ \
  x_{5}c_{11}=c_{11}\bar{x}_{5}. \]
  Using commutativity as well yields
\[c_{1}d x_{1}c_{3}\underline{r_{1}e_{1}e_{1}c_{4}f_1^{-1}}x_{2}t_{2,4}v_{2,4}c_7^{-1} c_{6}
 c_{7}\bar{x}_{3}f_3^{-1}\underline{r_{3} e_{3}e_{3}c_{8}f_3^{-1}}
x_{4}t_{4,6}v_{4,6}c_{11}^{-1}c_{10}
  c_{11}\bar{x}_{5}c_{13}^{-1}c_{12}c_{13}c_{13}c_{12}c_{13}c_{12}.\]
Following the same argument  given in \eqref{from-r-to-y} for the
underlined portions and renaming $c_7^{-1} c_{6}
 c_{7}=s_3,c_{11}^{-1}c_{10}c_{11}=s_5$ and $c_{13}^{-1}c_{12}c_{13}=s_6$
 we get
  \[c_{1}d
  x_{1}c_{3}y_2y_2c_4x_{2}t_{2,4}v_{2,4}s_3\bar{x}_{3}f_3^{-1}y_4y_4c_8
x_{4}t_{4,6}v_{4,6}s_5\bar{x}_{5}s_6c_{13}c_{12}c_{13}c_{12}.\]

\begin{figure}[htbp] \label{figure-alternate-lanterns-f2-f4-genus-6}
     \centering  \leavevmode
     \psfig{file=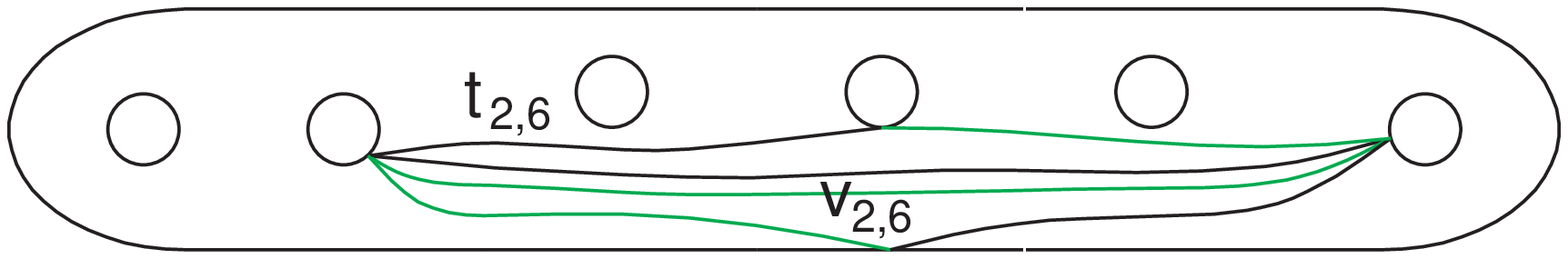,width=5.50in,clip=}
     \caption{}
 \end{figure}
 One more lantern substitution is needed
to eliminate $f_3$ and that is
\[f_3^{-1}=t_{2,6}v_{2,6}f_1^{-1}v_{2,4}^{-1}v_{4,6}^{-1}c_{13}^{-1}.\]
Result from substituting that is
 \[c_{1}d
  x_{1}c_{3}y_2y_2c_4x_{2}t_{2,4}v_{2,4}s_3\bar{x}_{3}t_{2,6}v_{2,6}f_1^{-1}v_{2,4}^{-1}v_{4,6}^{-1}c_{13}^{-1}y_4y_4c_8
x_{4}t_{4,6}v_{4,6}s_5\bar{x}_{5}s_6c_{13}c_{12}c_{13}c_{12}.\] The
idea in \eqref{x-conjugate} can be used for $t_{4,6}$ and $v_{4,6}$
 to write
\[t_{4,6}v_{4,6}=v_{4,6}\overline{t}_{4,6}.\]
Now, conjugating $y_4,c_8,x_4$ by $v_{4,6}$ and renaming them as
$y_{4,6}c_{8,12} x_{4,5}$, respectively, and likewise renaming
$v_{2,4}\bar{x}_{3}v_{2,4}^{-1}=\bar{x}_{3,2}$ and
$y_6=c_{13}^{-1}s_6c_{13}$ gives
\[c_{1}d
  x_{1}c_{3}y_2y_2c_4x_{2}t_{2,4}s_3\bar{x}_{3,2}t_{2,6}v_{2,6}f_1^{-1}y_{4,6}y_{4,6}c_{8,12}
x_{4,5}\overline{t}_{4,6}s_5\bar{x}_{5}y_6c_{12}c_{13}c_{12}.\] The
final lantern substitution is to replace $f_1^{-1}$ and it is
\[f_1^{-1}=t_{1,6}v_{1,6}c_1^{-1}c_3^{-1}v_{2,6}^{-1}c_{13}^{-1},\]
where $t_{1,6}$ and $v_{1,6}$ are defined in the same fashion.
Substituting that results in
\[c_{1}d
  x_{1}c_{3}y_2y_2c_4x_{2}t_{2,4}s_3\overline{x}_{3,2}t_{2,6}v_{2,6}t_{1,6}v_{1,6}c_1^{-1}c_3^{-1}v^{-1}_{2,6}c_{13}^{-1}y_{4,6}y_{4,6}c_{8,12}
x_{4,5}\overline{t}_{4,6}s_5\overline{x}_{5}y_6c_{12}c_{13}c_{12}.\]
After renaming two conjugations
$t_{1,2,6}=v_{2,6}t_{1,6}v^{-1}_{2,6}$ and
$w_6=c_{13}^{-1}y_6c_{13}$ following the braid relation
$c_{12}c_{13}c_{12}=c_{13}c_{12}c_{13}$ we will push $c_1^{-1}$ to
the right end in order to cancel it with the $c_1$ at the left end
of the next copy. Following the same idea for $c_3^{-1}$ gives
$u_1=c_3^{-1}dc_3$ after we invoke \eqref{x-conjugate} one more time
in order to write $x_{1}c_{3}=c_{3}\bar{x}_{1}.$ All of these
changes are realized in the final form of the genus 6 word that
follows:
\begin{equation}\label{genus-6-final}
\hspace{.5in}
\left(u_1\bar{x}_{1}y_2y_2c_4x_{2}t_{2,4}s_3\bar{x}_{3,2}t_{2,6}t_{1,2,6}v_{1,6}
  y_{4,6}y_{4,6}c_{8,12}
x_{4,5}\overline{t}_{4,6}s_5\bar{x}_{5}w_6c_{12}c_{13}\right)^3=1.
\end{equation}

\section{Applications} \label{applications}

In this section we will compute the homeomorphism  invariants of the
4-manifolds defined by the words in the previous section.  We will
denote by $X_g$ the manifolds that are given by the words
\eqref{rotation}, \eqref{genus-3-word-final}, and
\eqref{genus-4-rose-word-final} and those that are obtained from
them by  inserting  $k$  lantern relations will be denoted by $ X_{g,k}$.\\

\begin{prop} \label{signature-of-X} {  The signature and Euler characteristic of the Lefschetz
fibration $X_g,g=2,3,4,$ is given by
$\sigma(X_g)=-2\left(g+7\right)$ and $\chi \left(
X_{g}\right)=2g+22$, respectively.} \\
\end{prop}
\noindent \textbf{Proof:}
 \indent By checking the respective equations
we see that the number of cycles in those that define $X_2,X_3,$ and
$X_4$ is \ $3\left( 2g+6\right);$ therefore their Euler
characteristics are given by the formula
\begin{equation*}
 \chi \left( X_{g}\right) =2\left( 2-2g\right) +3\left(
2g+6\right) =\allowbreak 2g+22.
\end{equation*}
Here we used the well
known fact from the theory of Lefschetz fibrations  that the Euler
characteristic of
 a Lefschetz fibration $X^{4}\rightarrow S^{2}$ is given by the formula
\begin{equation} \label{euler-characteristic}
 \chi (X)=4-4g+s,
\end{equation}
 where $g$ is the genus of the fiber and $s$ is the number of
singular fibers, i.e., the number of vanishing cycles \cite{GS}.

 For
signature computations that  follow the reader is referred to
article \cite{ES}. First we compute $\sigma(X_2)$. \\ \indent Let
$C_2$ denote a chain of length 2 in $\mathcal{M}_2$, such as
$\left(c_{1}c_{2}\right)^6 \delta^{-1} $ and
$\left(c_{5}c_{4}\right)^6 \delta^{-1}$. Following the construction
of the word in \ref{genus-2-section} we have
\begin{eqnarray}\nonumber && C_2\cdot C_2^{-1}\\ \nonumber &=&\left(c_{1}c_{2}\right)^6 \delta^{-1}
\delta \left(c_{5}c_{4}\right)^{-6}=\left(c_{1}c_{2}\right)^6
 \left(c_{5}c_{4}\right)^{-6}\\
\nonumber &\equiv &\left( \left(c_{1}c_{2}\right)^2
 \left(c_{5}c_{4}\right)^{-2}\right)^{3}\ \ \ ( \text{commutativity} )\\
\nonumber &\equiv &\left( \left(c_{1}c_{2}\right)^2
  \delta^{-1} \left(c_{5}c_{4}\right)^6\left(c_{5}c_{4}\right)^{-2}\right)^{3}\ \ \ ( \text{chain relation } C_2 )\\
\nonumber &\equiv &\left( \left(c_{1}c_{2}\right)^2  \delta^{-1}
 \left(c_{5}c_{4}\right)^{4}\right)^{3}\ \ \ ( \text{cancelation } )\\
\nonumber &\equiv &\left( \left(c_{1}c_{2}\right)^2
xc_{3}c_{1}^{-2}c_{5}^{-2}
 \left(c_{5}c_{4}\right)^{4}\right)^{3}\ \ \ ( \text{lantern relation} )\\
 &\equiv & \cdots \cdots\\
\nonumber &\equiv & \left(c_{1}c_{2}
xc_{3}c_{4}c_{5}c_{5}c_{4}c_{5}c_{4}\right)^3\ \ \
(\text{commutativity, braid relations} )
 \end{eqnarray}

Cancelations do not change the signature and commutativity and braid
relations   have   zero signature (\cite{ES}, Proposition 3.6);
therefore we have
\begin{eqnarray*} \sigma \left( X_{2}\right) &=&I\left(
C_{2}\right) -I\left( C_{2}\right) +3I\left( C_{2}\right) +3I\left(
L\right) \\&=&-7-\left( -7\right) +3\left( -7\right) +3\left(
+1\right) \\&=&\allowbreak -18
\end{eqnarray*}

Next, we compute $\sigma(X_3).$ Let $C_2$ denote either of the two
chains $\left(c_{1}c_{2}\right)^6 \delta_{1}^{-1} $ or
$\left(c_{7}c_{6}\right)^6 \delta_{2}^{-1}$ of length 2 and $C_3$
denote the chain $\left(e_1c_{4}f_{1}\right)^4 \delta_{1}^{-1}
\delta_2^{-1}$ of length 3  in $\mathcal{M}_3$. Then construction of
the word in \ref{genus-3-section} gives

\begin{eqnarray}\nonumber && C_2\cdot C_3^{-1}\cdot C_2 \\ \nonumber &=&\left(c_{1}c_{2}\right)^6 \delta_{1}^{-1}
\delta_{1}\left(e_1c_{4}f_{1}\right)^{-4}  \delta_2\,
\delta_{2}^{-1}\left(c_{7}c_{6}\right)^{6} \\
\nonumber &\equiv &\left(
\left(c_{1}c_{2}\right)^2\left(e_1c_{4}f_{1}c_4\right)^{-1}
 \left(c_{7}c_{6}\right)^{2}\right)^{3}\ \ \ ( \text{commutativity and cancelations} )\\
\nonumber &\equiv &\left(
\left(c_{1}c_{2}\right)^2\left(e_1c_{4}f_{1}c_4\right)^{-1}\left(e_1c_{4}f_{1}c_4\right)^3
\delta_{1}^{-1} \delta_2^{-1}
 \left(c_{7}c_{6}\right)^{2}\right)^{3}\ \ \ ( \text{chain relation } C_3 )\\
\nonumber &\equiv &\left( \left(c_{1}c_{2}\right)^2\delta_{1}^{-1}
\left(e_1c_{4}f_{1}c_4\right)^2 \delta_2^{-1}
 \left(c_{7}c_{6}\right)^{2}\right)^{3}\ \ \ ( \text{commutativity and cancelations} )\\
\nonumber &\equiv &\left(
\left(c_{1}c_{2}\right)^2x_{1}c_{3}c_{1}^{-2}e_{1}^{-1}f_{1}^{-1}
\left(e_1c_{4}f_{1}c_4\right)^2
x_{2}c_{5}e_{1}^{-1}f_{1}^{-1}c_{7}^{-2}
 \left(c_{7}c_{6}\right)^{2}\right)^{3}\ \ \ ( \text{2 lantern relations }L )\\
\nonumber &\equiv &\cdots \cdots \\
\nonumber &\equiv & \left(c_{1}c_{2} x_{1}c_{3}rc_{8} c_{8}c_{4}
 x_{2}c_{5}c_{6}c_{7}\right)^3=1\ \ \
(\text{commutativity, braid relations} )
\end{eqnarray}
Keeping track of the relations that are used in the process we
obtain
\begin{eqnarray*} \sigma \left( X_{3}\right) &=&I\left(
C_{2}\right) -I\left( C_{3}\right)+I\left( C_{2}\right) +3I\left(
C_{3}\right) +3I\left( L\right)+3I\left( L\right) \\&=&-7-\left(
-6\right)+\left(-7\right) +3\left( -6\right) +3\left(
+1\right)+3\left( +1\right)
\\&=&\allowbreak -20
\end{eqnarray*}
Here we also used the fact that $\left( e_{1}c_{4}f_{1}\right)
^{4}=\left( e_{1}c_{4}f_{1}c_{4}\right) ^{3}$.

 We compute $\sigma
\left( X_{4}\right)$ last. Following its construction in
\ref{genus-4-section} we obtain

\begin{eqnarray}\nonumber && C_2\cdot C_2\cdot C_2\cdot E^{-1} \\ \nonumber &=&\left( b_{1}a_{1}\right) ^{6}\delta _{1}^{-1}
\left( b_{2}a_{2}\right) ^{6}\delta _{2}^{-1}\left(
b_{3}a_{3}\right) ^{6}\delta _{3}^{-1}\left( \alpha _{1}\alpha _{2}\alpha _{3}d\right) ^{-3}\delta _{1}\delta _{2}\delta _{3}\\
\nonumber &\equiv &\left( \left(b_{1}a_{1}\right)^2
\left(b_{2}a_{2}\right)^2 \left(b_{3}a_{3}\right)^2 \left( \alpha
_{1}\alpha _{2}\alpha _{3}d\right)^{-1}\right)^{3}\ \ \ ( \text{commutativity and cancelations} )\\
\nonumber &\equiv &\left(
 \left(b_{1}a_{1}\right)^2
\left(b_{2}a_{2}\right)^2 \left(b_{3}a_{3}\right)^2 \delta_{1}^{-1}
\delta_2^{-1} \delta_3^{-1} \left( \alpha _{1}\alpha _{2}\alpha
_{3}d\right)^{3} \left( \alpha
_{1}\alpha _{2}\alpha _{3}d\right)^{-1} \right)^{3}\ \ \ ( \text{star relation } E )\\
\nonumber &\equiv &\left(
 \left(b_{1}a_{1}\right)^2
\left(b_{2}a_{2}\right)^2 \left(b_{3}a_{3}\right)^2\delta_{1}^{-1}
\delta_2^{-1} \delta_3^{-1}  \left( \alpha
_{1}\alpha _{2}\alpha _{3}d\right)^{2}  \right)^{3}\ \ \ ( \text{commutativity and cancelations} )\\
\nonumber &\equiv &\left(
 \left(b_{1}a_{1}\right)^2
\left(b_{2}a_{2}\right)^2
\left(b_{3}a_{3}\right)^2x_{1}c_{1}a_{1}^{-1}a_{1}^{-1} \alpha
_{1}^{-1}\alpha _{2}^{-1}x_{2}c_{2}a_{2}^{-1}a_{2}^{-1} \alpha
_{2}^{-1}\alpha_{3}^{-1}\right. \\ \nonumber &&
\left.x_{3}c_{3}a_{3}^{-1}a_{3}^{-1}\alpha _{3}^{-1}\alpha _{1}^{-1}
\left( \alpha
_{1}\alpha _{2}\alpha _{3}d\right)^{2}  \right)^{3}\ \ \ ( \text{3 lantern relations }L )\\
\nonumber &\equiv &\cdots \cdots \\
\nonumber &\equiv & \left(a _{1}b _{1}a _{2}b _{2}a_{3}b
_{3}x_{1}c_{1}x_{2}c_{2} x_{3}c_{3}rd\right)^3=1\ \ \
(\text{commutativity, braid relations} )
\end{eqnarray}
From this we obtain
\begin{eqnarray*} \sigma \left( X_{4}\right) &=&3I\left(
C_{2}\right) -I\left( E\right)+3I\left(E\right) +3I\left(
L\right)+3I\left( L\right)+3I\left( L\right)
\\&=&3\left(-7\right)-\left( -5\right)+3\left(-5\right) +3\left(
+1\right) +3\left( +1\right)+3\left( +1\right)
\\&=&\allowbreak -22\\ && \hspace{4in} \square
\end{eqnarray*}

Consider now the fibrations $ X_{g,k}$  given by the words
\eqref{genus-2-lantern-inserted},\eqref{genus-3-alternate-word-final},
 and \eqref{genus-4-rose-word-modified} which are obtained from \eqref{rotation}, \eqref{genus-3-word-final},
and \eqref{genus-4-rose-word-final} by substituting $k$ lantern
relations.

\begin{prop} \label{lantern-blowup-prop} The Euler characteristic and the signature of the
manifold $X_{g,k}$ are given by
$\sigma\left(X_{g,k}\right)=\sigma\left(X_{g}\right)+k$ and
$\chi\left(X_{g,k}\right)=\chi\left(X_{g}\right)-k,g=2,3,4.$
\end{prop}

\noindent \textbf{Proof : } The only substitutions used in
\eqref{genus-2-lantern-inserted},\eqref
{genus-3-alternate-word-final},
 and \eqref{genus-4-rose-word-modified}that have
nonzero signature are lantern relations. The rest of the
modifications which result from commutativity and braid relations do
not have nonzero contributions (\cite{ES}, Proposition 3.6).
Cancelations also do not effect the signature. Since the signature
of each lantern relation is $+1$ half of the proof follows. The
other half follows from \eqref{euler-characteristic} and the fact
that each time we substitute a lantern relation the length of the
word reduces by one. \hfill $\square$

\begin{rem} \label{remark-on-k}
 To be more specific about $k$ we need to point out that $1\leq k\leq
6$ for genus $2$ and $1\leq k\leq 3$ for genus $3,4$. Therefore
$$-18\leq \sigma\left(X_{2,k}\right)\leq -12,\ \ -20\leq
\sigma\left(X_{3,k}\right)\leq -17, \ \ -22\leq
\sigma\left(X_{4,k}\right)\leq -19$$
\end{rem}
\begin{rem} In order to see that we have a positive relator for each
$k$ we will show what the word for $X_{2,1}$ becomes, for example.
$c_{1}c_{2} xc_{3}c_{4}c_{3}^{-2}khc_{1} c_{4}c_{5}c_{4}$ in the
third line of
 \eqref{genus-2-lantern-inserted} can be rewritten as $c_{1}c_{2}
xc_{4}^{-1}c_{3}c_{4}c_{3}^{-1}c_{4}c_{4}^{-1}kc_{4}c_{4}^{-1}hc_{4}c_{1}
c_{5}c_{4}$ and this becomes the positive relator $c_{1}c_{2}
xmnp\,c_{1} c_{5}c_{4},$ where
$m=c_{4}^{-1}c_{3}c_{4}c_{3}^{-1}c_{4},n=c_{4}^{-1}kc_{4},$ and
$p=c_{4}^{-1}hc_{4}$. Therefore the monodromy of  $X_{2,1}$ is
$$c_{1}c_{2} xmnp\,c_{1} c_{5}c_{4}\left(c_{1}c_{2}
xc_{3}c_{4}c_{5}c_{5}c_{4}c_{5}c_{4}\right)^2=1.$$ This is a
fibration with
$\sigma\left(X_{2,1}\right)=\sigma\left(X_{2}\right)+1=-17$ and
$\chi\left(X_{2,1}\right)=\chi\left(X_{2}\right)-1=25.$
\end{rem}

 \begin{rem} An interesting thing to observe here is the effect of
 substituting a lantern relation into the monodromy of $X_{g}$ on its
 homeomorphism invariants. Proposition \ref{lantern-blowup-prop} shows that it has
 the same effect on $X_g$ as that of a rational blowdown operation on it. Therefore it's
 an interesting  question to investigate whether or not   $X_{g}$ and $X_{g,k} \# k\overline{\mathbb{C}P}^{2}$
 are diffeomorphic.  See
 \cite{EG} for examples that answer this question in the negative.
\end{rem}
Next in our list is the word \eqref{genus-3-second-alternate-word}
 obtained from
  \eqref{genus-3-alternate-word-final} by substituting $m$ chain
  relations of length $3$ into $X_{3,k}$, which will be denoted by $X_{3,k,m},1 \leq m\leq k\leq 3$. This notation does not
  reflect the length of the chain for the sake of simplicity. Note that chain substitution must
follow a lantern substitution; therefore $m\leq k$.
  \begin{prop} $\sigma \left( X_{3,k,m}\right)    =-20+k-6m$ and $\chi\left(X_{g,k,m}\right)=28-k+10m$
  for $1\leq m\leq k \leq 3. $
  \end{prop}

  \noindent \textbf{Proof :} The signature of $X_{3}$ is $-20$ by
  Proposition  \eqref{signature-of-X}  and the signature of $X_{3,k}$ was found to
  be $-20+k$ in Proposition \eqref{lantern-blowup-prop}.
   Since $X_{3,k,m} $ is obtained from   $X_{3,k}$ by
substituting  $m$ chain relations of length $3$ and $C_3$ has
signature $-6$ (Proposition 3.10, \cite{ES}), we have
$$\sigma \left( X_{3,k,m}\right) =\sigma \left(
X_{3,k}\right)   +mI\left( C_{3}\right) =-20+k+m(-6),1\leq m\leq k
\leq 3.$$ Proposition  \eqref{signature-of-X}  gives
$\chi\left(X_{3}\right)=28$ and according to Proposition
\eqref{lantern-blowup-prop} $\chi\left(X_{g,k}\right)=28-k$. Since
 substitution of each $C_3$ results in increasing overall number of
cycles by $10$, its contribution to the Euler characteristic will be
$10$ according to \eqref{euler-characteristic}. Therefore we have
$\chi\left(X_{g,k,m}\right)=28-k+10m,1\leq m\leq k \leq 3. $

   \hspace{4.5in}
$\square$
\begin{rem} Possible values for $\sigma \left( X_{3,k,m}\right) $ are
$-23,-24,-25,-29,-30,-35$ and possible values for
$\chi\left(X_{g,k,m}\right)$ are $35,36,37,45,46,55$.
\end{rem}
Next, we will  compute the signatures of the achiral Lefschetz
fibrations \eqref{genus-even-achiral-final} and
\eqref{genus-odd-achiral-final}, denote them by $Z_g$. Assume,
first, $g$ is even and greater than $7$. $Z_g$ has  monodromy
\[c_{1}d
x_{1}c_{3}r_{1}e_{1}e_{1}c_{4}x_{2}c_{5}f_{2}^{-1} W_6 W_8\cdots W_g
c_{2g-2}
x_{g-1}c_{2g-1}c_{2g}c_{2g+1}c_{2g+1}c_{2g}c_{2g+1}c_{2g},\]

where  $W_i=c_{2i-6} x_{i-3}c_{2i-5}r_{i-3}e_{i-3}e_{i-3}c_{2i-4}
 x_{i-2}c_{2i-3}f_{i-2}^{-1},i=6,8,\dots,g$.

 From its construction
 in \ref{genus-g-even} we can see that this word originally contains two chains of
 length $2$, one on each end, and $g-2$ chains of length $3$ half of
 which are negatively oriented. Then we substituted $\displaystyle
 3(g-2)/2$ additional chains of length $3$ in order to replace
 the negatively oriented ones by positive exponents. We also substituted
 $3$ chains of length $2$ for the same reason. These substitutions resulted in
 $3(g-1)$ separating  negatively oriented boundary curves. Finally we
 introduced $3(g-1)$ lantern relations to eliminate them. The rest
 of the operations until we obtained
 \eqref{genus-even-achiral-final} are cancelations, commutativity
 and braid relations, which have zero contribution to the signature.
Combining all of that we can compute the signature of $Z_g$ as
$$\hspace{-.2in} \sigma \left( Z_{g}\right) =I\left( C_{2}\right) -\frac{g-2}{2}I\left( C_{3}\right)
+\frac{g-2}{2}I\left( C_{3}\right) -I\left(
C_{2}\right)+\frac{3\left( g-2\right) }{2}I\left( C_{3}\right)
+3I\left( C_{2}\right) +3\left( g-1\right) I\left( L\right)
$$
$$=0+\frac{3\left( g-2\right) }{2}(-6)+3(-7)+3(g-1)(+1)=-6g-6$$
Suppose now that $g$ is odd and greater than $6$. This time $Z_g$ is
given by the monodromy
\[c_{1}d
x_{1}c_{3}r_{1}c_{2g+2}c_{2g+2}c_{4}x_{2}c_{5}f_{2}^{-1}
 W_6 W_8\cdots W_{g-1} W_g,\]
where $W_i=c_{2i-6}
x_{i-3}c_{2i-5}r_{i-3}c_{2g+i-2}c_{2g+i-2}c_{2i-4}
 x_{i-2}c_{2i-3}f_{i-2}^{-1},i=6,8,\dots,g-1$ and \\
 $W_g=c_{2g-6} x_{g-3}c_{2g-5}r_{g-3}c_{3g-2}c_{3g-2}c_{2g-4}
x_{g-2}c_{2g-3}r_{g-2}f_{g-2}c_{3g-1}c_{2g-2} x_{g-1}c_{2g-1}
 c_{2g}c_{2g+1}$.
 Using a similar argument we calculate the signature of $Z_g$ as
 $$\sigma \left( Z_{g}\right) =I\left( C_{2}\right) -\frac{g-1}{2}I\left( C_{3}\right)
 +\frac{g-3}{2}I\left( C_{3}\right) +I\left( C_{2}\right)+\frac{3\left( g-1\right) }{2}I\left( C_{3}\right)
 +3\left( g-1\right) I\left( L\right) $$
 $$=-7-\frac{g-1}{2}\left( -6\right) +\frac{g-3}{2}\left( -6\right) +\left( -7\right)
 +\frac{3\left( g-1\right) }{2}\left( -6\right) +3\left(
 g-1\right)(+1)=-6g-2
 $$

 Note that $\sigma \left( Z_{g}\right)=\sigma \left( X_{g}\right)$
 for $g=2,3$. This is because the simplified form of the general
 construction leads to a positive relator. The existence of negative powers in the expression
 for higher genus, however, requires further substitution of lantern relations.
 We'll denote by $Y_{g}$  a genus $g$ Lefschetz fibration that
is obtained from either of the achiral Lefschetz fibrations
\eqref{genus-even-achiral-final} or \eqref{genus-odd-achiral-final}
by substituting into them a number of lantern relations until a
positive relator is obtained. In that regard the fibration given by
\eqref{final-word-genus-4} that is obtained from
\eqref{genus-4-word-with-negative} via $3$ lantern substitutions
will be denoted by $Y_{4}$. If $k$ additional lantern substitutions
are made into these positive words then the resulting manifold will
be denoted by $Y_{g,k}$. For example the positive relator
\eqref{genus-4-with-second-lantern-last} is denoted by $Y_{4,3}$
 because it is obtained via $3$ lantern substitutions into
 \eqref{final-word-genus-4}, which is equivalent to
 \eqref{genus-4-second-from-last}.
 We will now compute the  signatures  of $Y_{g},Y_{g,k} $.

\begin{prop}
$\sigma \left( Y_{4}\right)=-27,\sigma \left(
Y_{4,k}\right)=-27+k,\sigma \left( Y_{5}\right)=-29,\sigma \left(
Y_{6}\right)=-30.$
\end{prop}
\noindent \textbf{Proof :} $Y_{4}$ is obtained from $Z_{4}$ by
substituting lantern relation $3$ times; therefore
$$\sigma \left( Y_{4}\right)=\sigma \left( Z_{4}\right)+3(+1)=-6\cdot 4
-6+3=-27.$$ Similarly $ Y_{4,k} $ is obtained from $ Y_{4} $ by
substituting $k$ lantern relations,$1\leq k \leq 3$ ; therefore
$$\sigma \left( Y_{4,k}\right)=\sigma \left( Y_{4}\right)+k(+1)=-
27+k.$$ $Y_{5}$ is obtained from $Z_{5}$ by substituting lantern
relation $3$ times; therefore
$$\sigma \left( Y_{5}\right)=\sigma \left( Z_{5}\right)+3(+1)=-6\cdot
5 -2+3=-29.$$ A careful analysis shows that  $Y_{6}$ is obtained
from $Z_{6}$ by substituting lantern relation $12$ times ( $3$ for
each of the negative powers $f_2^{-1}, f_4^{-1}, f_3^{-1}, f_1^{-1}$
); therefore
$$\sigma \left( Y_{6}\right)=\sigma \left( Z_{6}\right)+12(+1)=-6\cdot
6-6+12=-30.$$ \hspace{5in} $\square$

We summarize what we found in the following table, which  includes
the Euler characteristic $\chi$, signature $\sigma$, the holomorphic
Euler characteristic $\displaystyle \chi _{h}=\frac{1}{4}\left(
\sigma +\chi
  \right)$, and the self-intersection of the first
Chern class $c_{1}^{2}=3\sigma +2\chi  $. The latter two are defined
for manifolds having almost complex structure and symplectic
Lefschetz fibrations are known to possess that.\\

\hspace{-.4in}
\begin{tabular}{l|c|c|c|c|c|c}
 & \ \  $\chi$ \ \ & \ \ $\sigma$ \  \ & \ $\chi _{h}$ \ \ & \ \ $c_{1}^{2}$ \ \ & $\pi_1$ &  \\
  \hline
$X_{2}$  & $26$ & $-18$ &  $2$ & $-2$ & 1 & \\   \hline
    $X_{2,k}$  & $26-k$ & $-18+k$ &  $2$ & $-2+k$ & 1 & \ \ $k=1,\dots,5$ \\
    \hline
      $X_{2,6}$  & $20$ & $-12$ &  $2$ & $4$ & \ \ $\mathbb{Z}_3$ &  \\
      \hline
  $X_{3}$  & $28$ & $-20$ & $2$ & $-4$  & 1 & \\  \hline
    $X_{3,k}$ & $28-k$ &  $-20+k$ & $2$ & $-4+k$ & 1 &\  $k=1,2,3$ \\
    \hline
     $X_{3,k,m}$ & $28-k+10m$ &  $-20+k-6m$ & $2+m$ & $-4+k+2m$ & 1 &\  $m,k=1,2,3,m\leq k$ \\
     \hline
      $X_{4}$  & $30$ & $-22$ & $2$ & $-6$  & 1 &\\  \hline
    $X_{4,k}$ & $30-k$ &  $-22+k$ & $2$ & $-6+k$ & 1 & \  $k=1,2,3$\\
    \hline
        $Y_{4}$  & $39$ & $-27$ & $3$ & $-3$ & 1 & \\  \hline
    $Y_{4,k}$ & $39-k$ &  $-27+k$ & $3$ & $-3+k$ & 1 & \  $k=1,2,3$\\
    \hline
        $Y_{5}$  & $41$ & $-29$ & $3$ & $-5$ & 1 & \\  \hline
        $Y_{6}$  & $46$ & $-30$ & $4$ & $2$  \\ \hline
\end{tabular}

\section*{Acknowledgment}
The authors thank Ronald J. Stern for helpful comments.   They also
thank Yusuke Kuno for the correction of the signature of the star
relation $E$ in Proposition 3.13 of \cite{ES} that was calculated as
$+5$.
\bibliographystyle{amsplain}

\addcontentsline{toc}{subsection}{BIBLIOGRAPHY}
\begin{thebibliography}{02}

\bibitem{ES}
H. Endo and S. Nagami, \textit{Signature of relations in mapping
class groups and non-holomorphic Lefschetz fibrations}, Trans. Amer.
Math. Soc. \textbf{357} (2005), 3179-3199.

\bibitem{EG}
H. Endo and Y.Gurtas, \textit{ Lantern Relations and Rational
Blowdowns}, arXiv:0808.0386.

\bibitem{Ge}
S. Gervais, \textit{Presentation and Central Extensions of Mapping
Class Groups}, Transactions of the American Mathematical Society,
Vol. \textbf{348} 8 (1996), 3097-3132.

\bibitem{GS}
R.Gompf and A.Stipsicz, \textit{An Introduction to 4-manifolds and
Kirby Calculus,} AMS Graduate Studies in Mathematics, {\bf 20}
(1999).




\bibitem{Oz}
B.Ozbagci, \textit{Signatures of Lefschetz fibrations}, Pacific
Journal of Mathematics, Vol. {\bf 202} 1 (2002), 99-118.

\bibitem{Lu3}
F. Luo, \textit{Torsion elements inthe mapping class group of a
surfaces}, preprint, arXiv:math.GT/0004048.



\end{thebibliography}

\end{document}